\documentclass[12pt]{amsart}

\usepackage{graphicx,amssymb,latexsym,amsfonts,txfonts,amsmath,amsthm}
\usepackage{pdfsync,color,tabularx,rotating}
\usepackage[all,cmtip]{xy}
\usepackage{url}
\usepackage{tikz}
\usepackage{amssymb}
\usepackage{comment}

\advance\hoffset by -0.9 truecm

\newtheorem{theo}{Theorem}[section]
\newtheorem{prop}[theo]{Proposition}

\newtheorem{lemm}[theo]{Lemma}

\theoremstyle{remark}

\theoremstyle{remark}

\theoremstyle{remark}

\theoremstyle{remark}

\begin{document}

\title{Regular maps with primitive automorphism groups}

\author{Gareth A. Jones and Martin Ma\v caj }

\address{School of Mathematical Sciences, University of Southampton, Southampton SO17 1BJ, UK}
\email{G.A.Jones@maths.soton.ac.uk}

\address{Faculty of Mathematics, Physics and Informatics, Comenius University, Bratislava, Slovakia}

\email{macaj@dcs.fmph.uniba.sk
}


% 05C10 topollogical aspects of graph theory, 20B15 primitive groups, 20B25 finite groups of automorphisms ..., 

\begin{abstract}
We classify the regular maps $\mathcal M$ which have automorphism groups $G$ acting faithfully and primitively on their vertices.
As a permutation group $G$ must be of almost simple or affine type, with dihedral point stabilisers.
We show that all such almost simple groups, namely all but a few groups ${\rm PSL}_2(q)$, ${\rm PGL}_2(q)$ and ${\rm Sz}(q)$,
arise from regular maps, which are always non-orientable. In the affine case, the maps $\mathcal M$ occur in orientable and non-orientable Petrie dual pairs.
We give the number of maps associated with each group, together with their genus and extended type.
Some of this builds on earlier work of the first author on generalised Paley maps, and on recent work of Jajcay, Li, \v Sir\'a\v n and Wang on maps with quasiprimitive automorphism groups. There are tables of data for the maps in appendices to this paper.
\end{abstract}

\maketitle

%%%%%%%%%%%%%%%%

\section{Introduction}\label{sec:Intro}

In~\cite{Jones} the first author classified the orientably regular maps $\mathcal M$ with orientation-preserving automorphism group ${\rm Aut}^+{\mathcal M}$ acting primitively and faithfully on their vertices. They are the generalised Paley maps, which are Cayley maps for the generalised Paley graphs introduced by Lim and Praeger in~\cite{LP}. The automorphism groups arising there are subgroups of ${\rm A\Gamma L}_1(q)$ for prime powers $q$, so they are all solvable, with easily described structure. Relaxing the condition of a faithful action simply yields regular cyclic coverings of these maps, branched  over the vertices, so the groups are again solvable.

Our aim here is to prove similar results for (fully) regular maps $\mathcal M$, with (full) automorphism group $G={\rm Aut}\,{\mathcal M}$ acting primitively on the vertices. The O'Nan--Scott Theorem divides finite primitive permutation groups into a small number of types, depending on their normal structure. In this particular  case the stabiliser of a vertex is a dihedral maximal subgroup $D$ of $G$. Now the only O'Nan--Scott types which allow dihedral point stabilisers are those consisting of the affine groups and the almost simple groups. We therefore first consider which primitive groups of these two types have dihedral point stabilisers. In the almost simple case a result of Li~\cite{Li} implies that only possibilities are the projective groups ${\rm PSL}_2(q)$ and ${\rm PGL}_2(q)$, for all except a few small values of $q$, and the Suzuki groups ${\rm Sz}(q)$ (see Theorem~\ref{th:LiThm}), while in the affine case this is equivalent to finding the irreducible representations of dihedral groups over prime fields (see Section~\ref{sec:affine}).

This action of $G$ is as the automorphism group of a regular map if and only if $G$ is generated by involutions $r_0$, $r_1$ and $r_2$ with $D=\langle r_1, r_2\rangle$ and $r_0r_2=r_2r_0$ (see Section~\ref{sec:regular}). All of these almost simple groups $G$ satisfy this condition, and the affine groups do provided $|G|$ is divisible by $4$, or equivalently the number and the valency of the vertices are not both odd. In the almost simple case all the maps are non-orientable (see Theorem~\ref{th:ASprimmaps}); in the affine case orientable and non-orientable maps both arise, in Petrie dual pairs, the former being generalised Paley maps or joins of chiral pairs of them (Theorem~\ref{th:affgps}).  As with orientably regular maps, relaxing the condition of a faithful action on the vertices yields only regular cyclic branched coverings of these maps (see Section~\ref{sec:non-faithful}). 

In all cases we give formulae which enumerate the regular maps $\mathcal M$ associated with a given group $G$ or, more precisely, with a given action of $G$, since some of these groups have two conjugacy classes of dihedral maximal subgroups. We explain in Section~\ref{sec:countingAS} how to determine the orientability of these maps, their extended type and their genus, giving examples in Sections~\ref{sec:countingAS} and \ref{sec:affexs} for almost simple and affine groups with small values of the relevant parameters.  In small cases, and for some infinite families, these properties can be obtained by hand, but in some larger cases we have used GAP~\cite{GAP} for this purpose. We have also verified our results for smaller groups and maps by comparing them with Conder's census of regular maps~\cite{Conder}, obtained using MAGMA and going up to genus $301$ and $602$ for orientable and non-orientable maps respectively.

Some of our results overlap with those of Jajcay, Li, \v Sir\'a\v n and Wang~\cite{JLSW} on regular and orientably regular maps with automorphism groups acting quasiprimitively on the vertices. Quasiprimitive groups form a much wider class than primitive groups, especially in the case of almost simple groups; in this case their results are necessarily only very partial, whereas assuming primitivity allows much more specific results to be obtained. In the case of affine groups our methods and results are similar to theirs, but we are able to give extra information about the maps, building on earlier work of the first author~\cite{Jones} to realise the orientable maps as generalised Paley maps or as joins of chiral pairs of such maps. Since  each of the affine groups arising is a semidirect product of a regular normal subgroup by the dihedral stabiliser $D$ of a vertex, each corresponding map $\mathcal M$ is a regular cover of a one-vertex regular map $\overline{\mathcal M}$ with automorphism group $D$; these quotient maps $\overline{\mathcal M}$ are described in Section~\ref{sec:dihedral}.

Appendices~A and B give tables of data about the maps with almost simple and affine automorphism groups respectively. A shorter version of this paper, without the Appendices, has been submitted for publication.

%A longer version of this paper is in preparation and will be made available online, containing extensive computer-generated tables of these maps.

%%%%%%%%%%%%%%%%%%%

\section{Regular maps}\label{sec:regular}

A group $G$ is isomorphic to the automorphism group ${\rm Aut}\,{\mathcal M}$ of a regular map $\mathcal M$ if and only if it has generators $r_i$ ($i=0, 1, 2$) satisfying
\begin{equation}\label{r_i}
r_i^2=1,\quad r_0r_2=r_2r_0.
\end{equation}
Given $\mathcal M$, its barycentric subdivision triangulates the underlying surface, and these generators $r_i$ can be realised as reflections in the sides of a triangle opposite a vertex or the midpoint of an edge or face. Conversely, given $G$ and $r_i$, $\mathcal M$ can be reconstructed by taking the vertices, edges and faces to be the cosets of the subgroups $\langle r_1, r_2\rangle$, $\langle r_0, r_2\rangle$ and $\langle r_0, r_1\rangle$, with incidence given by non-empty intersection.

A regular map $\mathcal M$ has empty boundary if and only if each $r_i\ne 1$, as we will assume without further comment. (The regular maps with non-empty boundary, all of them on the closed disc, are easily classified, see~\cite{JonesBdy} or~\cite{LS}; the only example with more then two vertices is the embedding of a cycle with $m$ vertices in the boundary of the disc, and this is vertex-primitive if and only if $m$ is prime.) Then $\mathcal M$ is orientable if and only if $G$ has a subgroup $G^+$ of index~$2$ containing no generator $r_i$.
The {\em type\/} of $\mathcal M$, in the notation of~\cite{CM}, is $\{m,n\}$ where $m$ and $n$ are the common valencies of its faces and vertices. Its {\em extended type\/} is $\{m,n\}_l$ where $l$ is the length of its Petrie polygons: these are closed zig-zag paths, turning alternately first left and first right at each vertex. By the Riemann--Hurwitz formula the Euler characteristic of $\mathcal M$ is
\[\chi=\frac{|G|}{2}\left(\frac{1}{n}-\frac{1}{2}+\frac{1}{m}\right),\]
and $\mathcal M$ has genus $g=1-\frac{\chi}{2}$ or $2-\chi$
as it is orientable or non-orientable.

Duality, transposing vertices and faces, is equivalent to transposing the generators $r_0$ and $r_2$, giving a regular map of type $\{n,m\}_l$ and the same orientability and genus. Petrie duality, transposing faces and Petrie polygons, is equivalent to transposing $r_0$ and $r_0r_2$, giving a regular map of type $\{l,n\}_m$; this preserves the embedded graph, but it can change the surface. Both operations preserve the automorphism group, and Petrie duality preserves its action on the vertices, including primitivity.

For simplicity, apart from a few comments in Section~\ref{sec:non-faithful}, we will assume that ${\rm Aut}\,\mathcal M$ acts faithfully on the vertices of $\mathcal M$, since regular maps with non-faithful actions on vertices, edges or faces have been thoroughly investigated by Li and \v Sir\'a\v n in~\cite{LS}. This assumption implies that $\mathcal M$ has valency $n>2$: the case $n=1$ is trivial, and if $n=2$ then $\mathcal M$ is an embedding of a cycle $C$ with $m$ vertices and edges, for some $m\ge 2$, so that by comparing orders one sees that the induced homomorphism ${\rm Aut}\,\mathcal M \to {\rm Aut}\,C\cong{\rm D}_m$ has a non-identity kernel. (In fact, the only regular maps of valency $n=2$ are the embedding $\{m,2\}$ of an $m$-cycle in the sphere and the antipodal quotient of $\{2m,2\}$ in the real projective plane, and these are vertex-primitive if and only if $m$ is prime.) 

%%%%%%%%%%%%%%%%%%%%DM}}

\section{Primitive permutation groups}\label{sec:primitive}

A permutation group is {\sl primitive\/} if the only equivalence relations it preserves are the identity and the universal relation. This is equivalent to the point-stabilisers being maximal subgroups.  The O'Nan--Scott Theorem (see~\cite[Ch.~4]{DM}) divides finite primitive permutation groups into a small number of classes, depending on their normal structure. The most important, here and in general, are those consisting of almost simple or affine groups. A group $G$ is {\em almost simple\/} if $T\le G\le {\rm Aut}\,T$ for some non-abelian finite simple group $T$, the socle of $G$; for comprehensive information about simple and almost simple groups see~\cite{ATLAS, Wilson}. A permutation group $G$ is of {\em affine type\/} if it has an elementary abelian normal subgroup $V$ acting regularly, in which case $V$ can be regarded as a $d$-dimensional vector space over ${\mathbb F}_p$ where $|V|=p^d$ with $p$ prime, and $G=V\rtimes G_0\le{\rm AGL}_d(p)$ with $G_0$, the stabiliser of $0\in V$, acting by conjugation on $V$ as a subgroup of ${\rm GL}_d(p)$; then $G$ is primitive if and only if $G_0$ acts irreducibly on $V$.

For any regular map $\mathcal M$, the stabiliser in $G:={\rm Aut}\,{\mathcal M}$ of a vertex of $\mathcal M$ must be a dihedral group of order $2n$, where $n$ is the valency of the vertices. We will restrict our attention to regular maps with more than two vertices, since transitive permutation groups of degree at most $2$ are trivially primitive. The following result shows that the only classes in the O'Nan--Scott Theorem which allow dihedral point-stabilisers are the almost simple and affine classes, so in this paper we will restrict our attention to them.

\begin{prop}
If a finite primitive permutation group $G$ has dihedral point-stabilisers $G_{\alpha}$ then $G$ is of affine or almost simple type.
\end{prop}

\noindent{\sl Proof.} (See~\cite[Ch.~4]{DM}), whose notation we follow, for the background to this proof.) By~\cite[Corollary 4.3B]{DM}, the socle $H$ of $G$ has the form $H=T_1\times \cdots \times T_m$ with $T_i\cong T$ for $i=1,\ldots, m$ where $T$ is a simple group. If $T$ is abelian then $H$ is elementary abelian and $G$ is of affine type, so we may assume that $T$ is nonabelian. If $H$ acts regularly, then by~\cite[Theorem 4.7(i)]{DM} $G_{\alpha}$ has no nontrivial solvable normal subgroups, so in particular it cannot be a dihedral group. If $H$ does not act regularly, then by~\cite[Theorem~4.6A]{DM} either $m=1$ and $G$ is almost simple, or $m\ge 2$ and $G$ is of diagonal or product type. In these last two cases either $G_{\alpha}$ contains a copy of $T$ (see \cite[Lemma~4.5B]{DM} and subsequent comments), which is impossible for a dihedral group, or $H_{\alpha}=R_1\times\cdots\times R_m$ where $1<R_i<T_i$ and each $R_i\cong R<T$ (see the proof of~\cite[Theorem~4.6A]{DM}). Since $H_{\alpha}$ is contained in a dihedral group $G_{\alpha}$ we must have $m=2$ and $|R|=2$, so $H_{\alpha}\cong {\rm V}_4$. Now $H$ is normal in $G$, so $H_{\alpha}$ is normal in $G_{\alpha}$, which is isomorphic to ${\rm D}_n$, and hence $n=2$ or $4$. If $n=2$ then $G_{\alpha}=H_{\alpha}$, so $G=H=T_1\times T_2$, which is imprimitive. Hence $n=4$, so $|G:H|=|G_{\alpha}:H_{\alpha}|=2$ and hence $G$ is isomorphic, as a permutation group, to $T\wr {\rm S}_2$ in its product action. However, this is also imprimitive, as no nonabelian simple group $T$ can have a maximal subgroup $R$ of order $2$. \hfill$\square$

%%%%%%%%%%%%%%%%%%%

\section{Almost simple automorphism groups}\label{sec:almost}

In this section we will consider automorphism groups of regular maps which act faithfully on the vertices as primitive permutation groups of almost simple type.
 Li~\cite[Lemma~3.1]{Li} has proved the following:

\begin{lemm}
Let $G$ be an almost simple primitive permutation group on $\Omega$ with socle $T$ (so that $G\le{\rm Aut}\,T$).  Assume that the point stabilisers $T_{\omega}\;(\omega\in\Omega)$ in $T$ are dihedral. Then $T$ is primitive, and one of the following holds:
\begin{itemize}
\item[(i)] $T \cong {\rm PSL}_2(q)$ and $T_{\omega}\cong {\rm D}_{(q\pm 1)/2}$, where $q$ is odd;
\item[(ii)] $T \cong {\rm PSL}(2, q)$ and $T_{\omega}\cong  {\rm D}_{q\pm 1}$, where $q$ is even;
\item[(iii)] $T \cong {\rm Sz}(q)$ and $T_{\omega} \cong {\rm D}_{q-1}$.
\end{itemize}
\end{lemm}

The automorphism groups of these groups $T$ are well known (see~\cite{ATLAS} or~\cite{Wilson}, for example). Using this, a straightforward inspection of the subgroups of ${\rm Aut}\,T$ containing $T$ gives a precise description of the almost simple primitive permutation groups $G$ with dihedral point stabilisers:

\begin{theo}\label{th:LiThm}
An almost simple primitive group $G$ has point stabilisers $D$ isomorphic to a dihedral group ${\rm D}_n$ if and only if one of the following occurs:
\begin{enumerate}
\item $G\cong{\rm PSL}_2(q)$ for odd $q>11$, with $n=(q-1)/2$;.
\item $G\cong{\rm PGL}_2(q)$ for odd $q>5,$ with $n=q-1$;
\item $G\cong{\rm PSL}_2(q)={\rm SL}_2(q)$ for $q=2^e\ge 4$, with $n=q-1$;
\item $G\cong{\rm PSL}_2(q)$ for odd $q>9$, with $n=(q+1)/2$;
\item $G\cong{\rm PGL}_2(q)$ for odd $q>5$, with $n=q+1$;
\item $G\cong{\rm PSL}_2(q)={\rm SL}_2(q)$ for $q=2^e\ge4$, with $n=q+1$;
\item $G\cong{\rm Sz}(q)$ for $q=2^e$ and odd $e>1$, with $n=q-1$ .
\end{enumerate}
In each case $G$ has a single conjugacy class of maximal subgroups $D\cong{\rm D}_n$ for the specified value of $n$.
\end{theo}

Note that while ${\rm PSL}_2(4)\cong{\rm PSL}_2(5)$, all groups $G$ listed here are mutually non-isomorphic.
They are all simple, except for the groups ${\rm PGL}_2(q)$ ($q$ odd) which contain ${\rm PSL}_2(q)$ with index $2$. In the excluded cases for small $q$, there are non-maximal subgroups $D\cong{\rm D}_n$: in case~(1), for instance, with $q=11$, each subgroup $D\cong {\rm D}_5$ is contained in a subgroup isomorphic to ${\rm A}_5$.

It is clear that if the automorphism group $G={\rm Aut}\,{\mathcal M}$ of a regular map $\mathcal M$ acts primitively and faithfully on its vertices, then $G$ and the vertex stabilisers $D$ must be as in one of the cases of Theorem~\ref{th:LiThm}. What is less clear (and at first seems rather surprising) is that all pairs $G$ and $D$ in Theorem~\ref{th:LiThm} arise in this way. This is less surprising if one recalls that the automorphism groups of regular maps are characterised by having generators $r_0$, $r_1$ and $r_2$ satisfying
\[r_i^2=(r_0r_2)^2=1.\]
We already have involutions $r_1$ and $r_2$ generating $D$, so by the maximality of $D$ it is sufficient to find an involution $r_0$ which is in the centraliser $C_G(r_2)$ of $r_2$ but not in $D$. In fact, by counting such involutions $r_0$ we will give a formula for the number $\nu_i(n)$ of isomorphism classes of $n$-valent regular maps $\mathcal M$ in each case (i) of Theorem~\ref{th:LiThm}, where $i=1,\ldots, 7$. By noting that $\nu_i(n)\ge 1$ for all relevant values of $n$, we deduce the following:

\begin{theo}\label{th:ASprimmaps}
Let $G$ be an almost simple group.
\begin{itemize}
\item[(a)] $G$ is the automorphism group of an $n$-valent regular map $\mathcal M$, acting primitively on its vertices, if and only if  $G$ and $n$ are as in Theorem~\ref{th:LiThm}.
\item[(b)] All such maps $\mathcal M$ are non-orientable.
\item[(c)] In each case $i=1,\ldots, 7$, the number $\nu_i(q)$ of such maps $\mathcal M$ is given by the following table.

\end{itemize}
\end{theo}
$$\begin{array}{|c|c|c|c|c|}\hline
\mbox{Case $i$} & \mbox{Group} & q\,{\rm mod}\,4 & n & \nu_i(q) \\ \hline
\mbox{(1)} & \mbox{PSL}_2\mbox{(}q\mbox{)} & 1 & (q-1)/2 & (q-5)\phi(n)/8e \\ \hline
\mbox{(1)} & \mbox{PSL}_2\mbox{(}q\mbox{)} & -1 & (q-1)/2 & (q+1)\phi(n)/8e \\ \hline
\mbox{(2)} & \mbox{PGL}_2\mbox{(}q\mbox{)} & \pm1 & q-1 & (q-2)\phi(n)/2e \\ \hline
\mbox{(3)} & \mbox{SL}_2\mbox{(}q\mbox{)} & 0 & q-1 & (q-2)\phi(n)/2e \\ \hline
\mbox{(4)} & \mbox{PSL}_2\mbox{(}q\mbox{)} & 1 & (q+1)/2 & (q-1)\phi(n)/8e \\ \hline
\mbox{(4)} & \mbox{PSL}_2\mbox{(}q\mbox{)} & -1 & (q+1)/2 & (q-3)\phi(n)/8e \\ \hline
\mbox{(5)} & \mbox{PGL}_2\mbox{(}q\mbox{)} & \pm1 & q+1 & (q-2)\phi(n)/2e \\ \hline
\mbox{(6)} & \mbox{SL}_2\mbox{(}q\mbox{)} & 0 & q+1 & (q-2)\phi(n)/2e \\ \hline
\mbox{(7)} & \mbox{Sz}\mbox{(}q\mbox{)} & 0 & q-1 & (q-2)\phi(n)/2e \\ \hline
\end{array}$$

\medskip

The next two sections are devoted to the proof of this result.

%%%%%%%%%%%%%%%%%%%%%%%%%%

\section{Counting maps}\label{sec:counting}

Let $G={\rm Aut}\,\mathcal M$, acting primitively and faithfully on the vertices of a regular map $\mathcal M$. Let $\mathcal O$ be an orbit of $A:={\rm Aut}\,G$ on maximal subgroups $D$ in $G$ isomorphic to ${\rm D}_n$. Our aim is to count the (isomorphism classes of) regular maps $\mathcal M$ with automorphism group $G$ and vertex-stabilisers $D\in{\mathcal O}$. We need slightly different arguments for odd and even values of $n$. We will then apply this general method of enumeration to the almost simple groups $G$ in cases (1) to (7) of Theorem~\ref{th:LiThm}.

\medskip

Suppose that $n$ is odd. Then each $D\in\mathcal O$ has a single conjugacy class of $n$ involutions, all non-central. The number of ordered pairs $(r_1, r_2)$ of these generating $D$ is $|{\rm Aut}\,D|=n\phi(n)$. If $N:=N_A(D)$ and $C:=C_A(D)$ are the normaliser and centraliser of $D$ in $A$ then $N/C$ acts semiregularly on these generating pairs, with $n\phi(n)/|N:C|=n\phi(n)|C|/|N|$ orbits of length $|N:C|$. For such pairs the number $I$ of involutions in $C_G(r_2)$ is independent of $r_2$; only one of these involutions (namely $r_2$) is in $D$, so there are $I-1$ involutions $r_0\in G\setminus D$ commuting with $r_2$, and each of these, together with $r_1$ and $r_2$, generates $G$ and hence determines a regular map $\mathcal M$ with ${\rm Aut}\,\mathcal M\cong G$. Since $C$ acts semiregularly on such elements $r_0$, it has $(I-1)/|C|$ orbits of length $C$ on them. This applies to each of the $n\phi(n)|C|/|N|$ orbits of $N$ on generating pairs $(r_1, r_2)$ for $D$, so the number of isomorphism classes of maps $\mathcal M$ arising from the orbit $\mathcal O$ is
\[\frac{n\phi(n)|C|}{|N|} \cdot \frac{(I-1)}{|C| }= \frac{n\phi(n)(I-1)}{|N|}.\]

\medskip

Suppose now that $n$ is even. Since we are assuming that $n>2$, each $D\in\mathcal O$ has a unique central involution $z:=(r_1r_2)^{n/2}$, and $D=C_G(z)$ by the maximality of $D$. In addition to $z$ there are $n$ non-central involutions in $D$, and the number of ordered pairs $(r_1, r_2)$ of these generating $D$ is $|{\rm Aut}\,D|=n\phi(n)$. If $N:=N_A(D)$ and $C:=C_A(D)$ are the setwise and pointwise stabilisers of $D$ in $A$ then $N/C$ acts semiregularly on these generating pairs, with $n\phi(n)/|N:C|=n\phi(n)|C|/|N|$ orbits of length $|N:C|$. Let $I$ denote the common number of involutions in $C_G(r_2)$ for such elements $r_2$; then three of these involutions (namely $r_2$, $z$ and $r_2z$) are in $D$, so $I-3$ are outside $D$ and yield maps $\mathcal M$ with automorphism group $G$. As when $n$ is odd, $C$ acts semiregularly on these elements $r_0$, so arguing as in Case~A we see that the number of maps obtained from $\mathcal O$ is
\[\frac{n\phi(n)|C|}{|N|} \cdot \frac{(I-3)}{|C| }= \frac{n\phi(n)(I-3)}{|N|}.\]

Now suppose that the $n$ non-central involutions $r_2\in D$ lie in distinct orbits under $A$, necessarily two orbits, with $n/2$ in each, so there are $n\phi(n)/2$ ordered generating pairs $r_1, r_2$ for $D$ with $r_2$ in each orbit.
Let $I^+$ and $I^-$ be the numbers of involutions in $C_G(r_2)$ for involutions $r_2$ in these orbits. In each case three of these involutions (namely $r_2$, $z$ and $r_2z$) are in $D$, so $I^{\pm}-3$ of them are outside $D$. Since $C$ acts semiregularly on such elements $r_0\in G\setminus D$, it has $(I-3)/|C|$ orbits of length $C$ on them. This applies to each of the $n\phi(n)|C|/2|N|$ orbits of $N$ on generating pairs $(r_1, r_2)$ for $D$, so the number of maps arising from the orbit $\mathcal O$ is
\[\frac{n\phi(n)(I^+-3)}{2|N|}+\frac{n\phi(n)(I^--3)}{2|N|}=\frac{n\phi(n)(I^++I^--6)}{2|N|}.\]

%%%%%%%%%%%%%%%%%%%%%%%%%%%

\section{Counting maps for almost simple groups}\label{sec:countingAS}

We will now apply these above general arguments to the various families of almost simple groups, dealing first with cases where $n$ is odd, and starting with one of the easiest of those. First we need some basic facts about these groups $G$. 

In $G={\rm PSL}_2(q)$, where $q$ is an odd prime power, for each $n=(q\pm 1)/2$ there is a single conjugacy class of subgroups $D\cong{\rm D}_n$ in $G$, and they are maximal for all $q\ge 11$, except for subgroups isomorphic to ${\rm D}_5$ in ${\rm PSL}_2(11)$. The normaliser of $D$ in ${\rm PGL}_2(q)$, which contains $G$ with index $2$, is a dihedral group $D^*$ of order $4n$. Extending $D^*$ by the Galois group ${\rm C}_e$ of ${\mathbb F}_q$, where $q=p^e$ for some prime $p$, gives the normaliser $N=N_A(D)$ of $D$ in $A:={\rm Aut}\,G={\rm P\Gamma L}_2(q)$, so $|N|=4ne$. All involutions in $G$ are conjugate, with centralisers isomorphic to ${\rm D}_{(q\pm 1)/2}$, whichever has order divisible by $4$; this contains
$I=\frac{1}{2}(q\pm 1)+1$
involutions.

In $G={\rm PGL}_2(q)$, where $q$ is an odd prime power, for each $n=q\pm 1$ there is a single conjugacy class of subgroups $D\cong{\rm D}_n$ in $G$, and they are maximal for all $q>5$. There are two conjugacy classes of involutions in $G$, one contained in $G^+:={\rm PSL}_2(q)$ and the other in $G^-:=G\setminus G^+$.  Involutions $i\in G^+$ have centralisers $C_G(i)\cong {\rm D}_{q-1}$ or ${\rm D}_{q+1}$ as $q\equiv 1$ or $-1$ mod~$(4)$, that is, $C_G(i)\cong{\rm D}_k$ where $k=q\pm 1\equiv 0$ mod~$(4)$, whereas involutions $i\in G^-$ have centralisers $C_G(i)\cong {\rm D}_{q+1}$ or ${\rm D}_{q-1}$ as $q\equiv 1$ or $-1$ mod~$(4)$, that is, $C_G(i)\cong{\rm D}_l$ where $l=q\pm 1\equiv 2$ mod~$(4)$. Since $k$ and $l$ are both even, $C_G(i)$ contains $I^+=k+1$ or $I^-=l+1$ involutions as $i\in G^+$ or not, with three of them in $D$, so $I^{\pm}-3=k-2$ or $l-2$ of them are outside $D$. We have $|N|=2ne$.

In $G={\rm SL}_2(q)$ for $q=2^e\ge 4$, there is a single conjugacy class of subgroups $D\cong{\rm D}_n$ for each $n=q\pm 1$, and they are maximal for all $q\ge 4$. All involutions in $G$ are conjugate, with centraliser an elementary abelian $2$-group ${\rm V}_q$ of order $q$ containing $I=q-1$ involutions. We have $|N|=2ne$.

In $G={\rm Sz}(q)$ (see~\cite[\S 4.2]{Wilson} for these groups) for $q=2^e$ (odd $e\ge 3$), with $n=q-1$, there is a single conjugacy class of subgroups $D\cong{\rm D}_n$, and they are maximal for all such $q$. All involutions in $G$ are conjugate, with centraliser a Sylow $2$-subgroup of order $q^2$ containing $I=n=q-1$ involutions. {We have $|N|=2ne$.

We will now use this information to give a formula (or formulae) for the number $\nu_i(q)$ of maps $\mathcal M$ associated with a group $G=G(q)$ in each case $(i)$. We will briefly describe the maps for a few small values of $q$, and where possible we will identify them in Conder's census of regular maps~\cite{Conder}.

Before doing this we will explain how to determine the orientability and type, and hence the genus, of each map, allowing such descriptions and identifications to be made. Using explicit standard generators $r_i$ for $G$,
one can calculate the face-valency $m$ and the Petrie length $l$ of $\mathcal M$ as the orders of the elements $r_0r_1$ and $r_0r_1r_2$. (In our cases, the vertex-valency is known, as the chosen value of $n$.)}
In cases (1) to (6) this can be done by hand, but as one might expect, the Suzuki groups ${\rm Sz}(q)$ are technically much harder to consider: they are represented by $4\times 4$ matrices, rather than $2\times2$, and  since the smallest Suzuki group has order $2912$0 the corresponding genera are well outside the range of~\cite{Conder} (the smallest is $2290$).
One can obtain canonical generators for all regular maps with
automorphism group ${\rm PSL}_2(q)$, ${\rm PGL}_2(q)$ or ${\rm Sz}(q)$ from~\cite{CPS} or~\cite{DJ}. Moreover,
in the case of Suzuki groups and in cases~(1) to (6) when $q$ is large, one can use
computational tools such as GAP~\cite{GAP} in order to obtain information.
%{\color{blue}[In the ADAM version, this paragraph replaced Section~\ref{sec:type}. For the arXiv version, retain this paragraph or include Section~\ref{sec:type}?]}

\medskip

{\bf Case (1)} Let $G={\rm PSL}_2(q)$, with $n=(q-1)/2$. If $q\equiv 1$ mod~$(4)$, so that $n$ is even, the number of regular maps $\mathcal M$ arising is
\[\nu_1(q)=\frac{n\phi(n)(I-3)}{|N|} = \frac{(n-2)\phi(n)}{4e}=\frac{q-5}{8e}\phi\left(\frac{q-1}{2}\right).\]
For instance, if $q=13$, so that $e=1$, we obtain two maps; these appear in Conder's census~\cite{Conder} as the Petrie dual pair N106.11 and N142.11 of types $\{7,6\}_{13}$ and $\{13,6\}_7$.
If $q=17$, so that $n=8$, we obtain six maps: these are
\begin{itemize}
\item  N325.6 with Petrie dual N325.7 and self-Petrie dual map N325.8 all of type $\{9,8\}_9$;
\item  N325.9 of type $\{9,8\}_{17}$, and its Petrie dual N389.5 of type $\{17, 8\}_9$;
\item a self-Petrie dual map  N389.6 of type $\{17, 8\}_{17}$.
\end{itemize}

If $q\equiv -1$ mod~$(4)$, so that $n$ is odd, a similar argument gives
\[\nu_1(q)=\frac{n\phi(n)(I-1)}{|N|} = \frac{(n+1)\phi(n)}{4e}=\frac{q+1}{8e}\phi\left(\frac{q-1}{2}\right).\]
Thus if $q=19$ we obtain $15$ maps $\mathcal M$. These are
\begin{itemize}
\item N97.1 of type $\{3,9\}_9$ and its Petrie dual N477.6 of type $\{9,9\}_3$, which is self-dual;
\item N97.2 of type $\{3,9\}_{10}$ and its Petrie dual N496.3 of type $\{10,9\}_3$;
\item a self-Petrie-dual map N325.3 of type $\{5,9\}_5$;
\item N325.4 of type $\{5,9\}_{10}$ and its Petrie dual N496.4 of type $\{10,9\}_5$;
\item a dual pair N477.8 (one of them is self-Petrie-dual) and a self-dual map N477.7  (whose Petrie-dual is the second map N477.8), all of type $\{9,9\}_9$
\item a dual pair N477.9 of type $\{9,9\}_{10}$ and their Petrie duals N496.5 and N496.6
of type $\{10,9\}_9$;
\item a self-Petrie-dual map N496.7 of type $\{10,9\}_{10}$.
\end{itemize}

\medskip

{\bf Case (2)} Let $G={\rm PGL}_2(q)$, with $n=q-1$ where $q$ is odd. The total number of maps $\mathcal M$, with $r_2\in G^+$ or $G\setminus G^+$, is
\[\nu_2(q)=\frac{n\phi(n)(k-2)}{4ne}+\frac{n\phi(n)(l-2)}{4ne}=\frac{\phi(n)(k+l-4)}{4e}=\frac{(q-2)\phi(q-1)}{2e}.\]
since $(k+l)/2=q$.
For instance, when $q=7$ this gives five maps, of which three have $r_2\in G^+$ and two have $r_2\in G^-$. From~\cite{Conder} we find the following:
\begin{itemize}
\item a Petrie dual pair N16.2 of type $\{4,6\}_8$ and N37.3 of type $\{8,6\}_4$;
\item a Petrie dual pair N30.5 of type $\{6,6\}_8$ and N37.4 of type $\{8,6\}_6$;
\item a self-Petrie-dual map N34.5 of type $\{7,6\}_7$. 
\end{itemize}
When $q=9$ we have $7$ maps; they are
\begin{itemize}
\item a Petrie dual pair N17.1 of type $\{3,8\}_{10}$ and N101.9 of type $\{10,8\}_3$;
\item a Petrie dual pair N47.2 of type $\{4,8\}_5$ and N65.3 of type $\{5,8\}_4$;
\item a Petrie dual pair N65.2 of type $\{5,8\}_8$ and N92.5 of type $\{8,8\}_5$ (self-dual);
\item a self-Petrie-dual map N101.8 of type $\{10,8\}_{10}$. 
\end{itemize}

\medskip

{\bf Case (3)} Let $G={\rm SL}_2(q)$ for $q=2^e\ge 4$, with $n=q-1$. Then
\[\nu_3(q)=\frac{n\phi(n)(I-1)}{|N|} = \frac{(n-1)\phi(n)}{2e} = \frac{(q-2)\phi(q-1)}{2e}.\]
If $q=4$ there is one map, the antipodal quotient of the dodecahedron, of type $\{5,3\}_5$,
an embedding of Petersen's graph in the real projective plane.
If $q=8$ we obtain six maps: 
\begin{itemize}
\item a Petrie dual pair N8.1 of type $\{3,7\}_9$ and N64.3 of type $\{9,7\}_3$;
\item a dual pair N56.5 of type $\{7,7\}_9$ and their Petrie duals N64.4 and N64.5 of type $\{9,7\}_7$.
\end{itemize}

If $q=16$ or $32$ we obtain $14$ or $90$ maps.

\medskip

{\bf Case (4)} Let $G={\rm PSL}_2(q)$ with $n=(q+1)/2$. If $q\equiv 1$ mod~$(4)$, so that $n$ is odd, we have
\[\nu_4(q)=\frac{n\phi(n)(I-1)}{|N|} = \frac{(n-1)\phi(n)}{4e}=\frac{q-1}{8e}\phi\left(\frac{q+1}{2}\right).\]
For instance, if $q=13$ we obtain nine maps $\mathcal M$. They appear in~\cite{Conder} as the following:
\begin{itemize}
\item N15.1 of type $\{3,7\}_{13}$ and its Petrie dual N155.5 of type $\{13,7\}_3$;
\item N106.11 of type $\{6,7\}_{13}$  and its Petrie dual N155.7 of type $\{13,7\}_6$;
\item N155.6 of type $\{13,7\}_{13}$ (self-Petrie-dual);
\item N119.7 with trinity symmetry (that is, both self-dual and self-Petrie-dual), together with a dual pair N119.5 (one of them self-Petrie-dual) and a self-dual map N119.6  (whose Petrie-dual is N119.5), all of type $\{7,7\}_7$.
\end{itemize}

If $q=17$ there are $12$ maps.
If $q\equiv -1$ mod~$(4)$ we have
\[\nu_4(q)=\frac{n\phi(n)(I-3)}{|N|} = \frac{(n-2)\phi(n)}{4e}=\frac{q-3}{8e}\phi\left(\frac{q+1}{2}\right).\]
For instance, if $q=11$ we obtain two maps; these are the Petrie dual pair  N46.6 and N57.4 of types $\{5,6\}_6$ and $\{6,6\}_5$ (the latter self-dual).
If $q=19$ we obtain eight maps.

\medskip

{\bf Case (5)} Let $G={\rm PGL}_2(q)$, with $n=q+1$ where $q$ is odd. The calculation is identical to that for $n=q-1$ in case~(2), except that at the end  we have
\[\nu_5(q)=\frac{\phi(n)(k+l-4)}{4e}=\frac{(q-2)\phi(q+1)}{2e}.\]
For instance, when $q=7$ this gives ten maps, of which six have $r_2\in G^+$ and four have $r_2\in G^-$.
From~\cite{Conder} we find the following:
\begin{itemize}
\item a Petrie dual pair N9.1 of type $\{3,8\}_7$ and N41.3 of type $\{7,8\}_3$,
\item a Petrie dual pair N9.2 of type $\{3,8\}_8$ and N44.5 of type $\{8,8\}_3$ (self-dual),
\item a Petrie dual pair N23.2 of type $\{4,8\}_7$ and N41.2 of type $\{7,8\}_4$,
\item a Petrie dual pair N23.3 of type $\{4,8\}_6$ and N37.3  of type $\{6,8\}_4$,
\item self-Petrie-dual maps N37.4 of type $\{6,8\}_6$ and N44.6 of type $\{8,8\}_8$ (self-dual).
\end{itemize}
The first map N9.1 can be obtained from Klein's Hurwitz map R3.1 of type $\{3,7\}_8$, on his quartic curve of genus $3$, by applying Wilson's`opposite' operation $DPD=PDP$, cutting the map along its edges and then rejoining the edges with the reverse orientation. Similarly, N41.3 can be obtained by applying the triality operation $PD$ to R3.1.

\medskip

{\bf Case (6)} Let $G={\rm SL}_2(q)$ for $q=2^e\ge 4$, with $n=q+1$. Then
\[\nu_6(q)=\frac{n\phi(n)(I-1)}{|N|} = \frac{(n-3)\phi(n)}{2e} = \frac{(q-2)\phi(q+1)}{2e}.\]
If $q=4$, so that $G\cong{\rm A}_5$, we obtain two  maps: the antipodal quotient of the icosahedron, of genus~$1$ and type $\{3,5\}_5$, and its Petrie dual N5.3, the antipodal quotient of the great dodecahedron, of type $\{5,5\}_3$ and genus $5$. If $q=8$ we obtain six maps:
\begin{itemize}
\item N16.1 of type $\{3,9\}_7$ and its Petrie dual N64.3 of type $\{7,9\}_3$;
\item N64.4 and its Petrie dual N64.5 of type $\{7,9\}_7$;
\item a dual pair N72.9 of type $\{9,9\}_9$, which is also a Petrie dual pair.
\end{itemize}
If $q=16$ or $32$ we obtain $28$ or $60$ maps.

\medskip

{\bf Case (7)} Let $G={\rm Sz}(q)$ for $q=2^e$ (odd $e\ge 3$), with $n=q-1$. Then
\[\nu_7(q)=\frac{n\phi(n)(I-1)}{|N|} = \frac{(n-1)\phi(n)}{2e}  = \frac{(q-2)\phi(q-1)}{2e},\]
the same as for ${\rm PGL}_2(q)$ and ${\rm SL}_2(q)$ with $n=q-1$. If $q=8$ we obtain six maps, all of genus too large to be included in~\cite{Conder}.
However, a computer search shows that there are a map of type $\{5, 7\}_{13}$ and genus
$2290$ together with its Petrie dual of type $\{13,7\}_5$ and genus $4082$, a Petrie dual and
simultaneously dual pair of maps of type $\{7, 7\}_7$ and genus $3122$, and a Petrie dual pair
of maps $\{13, 7\}_{13}$ and genus $4082$. (In fact there are 14 regular maps for Sz(8), shown in the following table, each
obtained from these by using duality and Petrie duality; we list these maps, but not those corresponding to the almost simple groups considered earlier, 
because of the extra complexity and possible unfamiliarity of the Suzuki groups.)

$$\begin{array}{|c|c|c|}\hline
\hbox{Type} & \hbox{Genus} & \hbox{Number of maps} \\ \hline
 \{  5,  7 \}_{13} & 2290 & 1 \\ \hline
 \{  5, 13 \}_{ 7} & 3250 & 1 \\ \hline
 \{  7,  5 \}_{13} & 2290 & 1 \\ \hline
 \{  7,  7 \}_{ 7} & 3122 & 2 \\ \hline
 \{  7, 13 \}_{ 5} & 4082 & 1 \\ \hline
 \{  7, 13 \}_{13} & 4082 & 2 \\ \hline
 \{ 13,  5 \}_{ 7} & 3250 & 1 \\ \hline
 \{ 13,  7 \}_{ 5} & 4082& 1 \\ \hline
 \{ 13,  7 \}_{13} & 4082 & 2 \\ \hline
 \{ 13, 13 \}_{ 7} & 5042 & 2 \\ \hline
\end{array}$$
  
 \medskip

 When $q=32$ the above formula gives $90$ maps, out of a total of $186$ regular maps for ${\rm Sz}(32)$, again confirmed by a computer search.

\medskip

In each of the cases $i=1,\ldots, 7$ we have $\nu_i(q)\ge 1$ for all relevant values of $q$, so this completes the proof of Theorem~\ref{th:ASprimmaps}(a) and (c). For Theorem~\ref{th:ASprimmaps}(b) we use the fact that $\mathcal M$ is orientable if and only if $G$ has a subgroup $G^+$ of index $2$ with each $r_i\not\in G^+$. In the case of the simple groups ${\rm PSL}_2(q)$ and ${\rm Sz}(q)$ there is clearly no such subgroup. When $G={\rm PGL}_2(q)$ there is a unique subgroup $G^+$ of index $2$, namely ${\rm PSL}_2(q)$; however, since $\langle r_1r_2\rangle = D^+\not\le G^+$ in cases~(2) and (5), the generators $r_1$ and $r_2$ must lie in different cosets of $G^+$, so again the maps are all non-orientable. This completes the proof of Theorem~\ref{th:ASprimmaps}.

\section{Affine groups and their associated maps}\label{sec:affine}

Vertex-primitive maps $\mathcal M$ with affine automorphism groups $G$ can be enumerated and described by using the methods developed in Section~\ref{sec:counting}, or by obtaining the required information from~\cite{JLSW}. However, it is more advantageous for us to use the results in~\cite{Jones}, since these allow us to identify the orientable maps $\mathcal M$ as regular generalised Paley maps, or as joins of chiral pairs of generalised Paley maps; the non-orientable maps $\mathcal M$ are then obtained as their Petrie duals.

If $G$ is an affine group then since it acts primitively and faithfully on the vertices of $\mathcal M$, the vertex set $V$ can be identified with the socle of $G$: this is an elementary abelian $p$-group of order $q=p^d$ for some prime $p$, or equivalently a $d$-dimensional vector space over ${\mathbb F}_p$. Then $G=V\rtimes D$ with $V$ acting regularly by translations, and the stabiliser $G_0=D=\langle a, b\mid a^n=b^2=(ab)^2=1\rangle\cong{\rm D}_n$ of the vertex $0\in V$ acts faithfully on $V$ as an irreducible subgroup of ${\rm GL}_d(p)$. (Recall that we assume that $|V|>2$.) Since we are assuming that $n>2$, there is a unique cyclic subgroup $D^+=\langle a\rangle\cong{\rm C}_n$ of index~$2$ in $D$, preserving the local orientation around $0$. Then $G$ has a subgroup $V\rtimes D^+$ of index $2$; it is unique if $n$ is odd, but if $n$ is even there are two others, isomorphic to $V\rtimes{\rm D}_{n/2}$. Let us define $D^-:=D\setminus D^+$ and $G^-:=G\setminus G^+$, the nontrivial cosets of $D^+$ in $D$ and of $G^+$ in $G$.

If $r_1$ and $r_2$ are the standard involutions generating $D$, define $I(r_2)$ to be the set of involutions in $C_G(r_2)\setminus D$, and let $I^{\pm}(r_2)=I(r_2)\cap G^{\pm}$. By the maximality of $D$, if $r_0\in I(r_2)$ then $G=\langle r_0,r_1, r_2\rangle$. Since $r_0^2=1$ and $r_0r_2=r_2r_0$, this triple ${\bf r}=(r_i)$ determines a regular map ${\mathcal M}={\mathcal M}({\bf r})$ with ${\rm Aut}\,{\mathcal M}\cong G$ acting primitively on $V$. If $r_0\in I(r_2)$ then $r_0r_2\in I(r_2)$ also, so the triple ${\bf r}'=(r_0r_2, r_1, r_2)$ determines a regular map ${\mathcal M}'={\mathcal M}({\bf r}')$, the Petrie dual $P({\mathcal M})$ of ${\mathcal M}$. 

We have $r_1, r_2\in D^-\subset G^-$, whereas if $n$ is even they lie in different cosets of the other two subgroups of index $2$ in $G$. Thus $\mathcal M$ is orientable if and only if $r_0\in G^-$. Now $r_0$ and $r_0r_2$ lie in different cosets of $G^+$, so by choosing $r_0\in I^-$ we obtain an orientable map $\mathcal M$ which we will denote by $\mathcal M^+$, with a non-orientable Petrie dual $\mathcal M'$, which we will denote by $\mathcal M^-$, corresponding to $r_0r_2\in I^+(r_2)$. Thus the affine maps arise in Petrie dual pairs $\mathcal M^{\pm}$, with $\mathcal M^+$ orientable and $\mathcal M^-$ non-orientable, so it is sufficient to classify and describe the maps $\mathcal M^+$, since the remaining maps $\mathcal M^-$ are easily found as their Petrie duals.

Each map $\mathcal M^+$ is orientably regular, with orientation-preserving automorphism group $G^+=V\rtimes D^+$. Since $D$ acts irreducibly on $V$, and $D^+$ is a normal subgroup of index $2$ in $D$, there are just two possibilities for the action of $D^+$ on $V$:
\begin{itemize}
\item[(A)] $V$ is an irreducible module for $D^+$, or
\item[(B)] $V=V_0\oplus V_1$ where $V_0$ and $V_1$ are irreducible $D^+$-submodules transposed by the elements of $D\setminus D^+$, with mutually inverse actions of $D^+$.
\end{itemize}

In case~A, $G^+$ acts primitively and faithfully on $V$. Now the orientably regular maps $\mathcal M$ with orientation-preserving automorphism groups ${\rm Aut}^+\mathcal M$ acting primitively and faithfully on their vertices have been classified in~\cite{Jones} as generalised Paley maps (which we will define shortly), so it sufficient to determine which of these are fully regular. In case~B, $\mathcal M^+$ is the join (or minimal regular cover) of the chiral pair of orientably regular maps $\mathcal M_i=\mathcal M^+/V_i\;(i=0, 1)$; by the irreducibility of the quotient modules $V/V_i\cong V_{1-i}$, each $\mathcal M_i$ has orientation-preserving automorphism group ${\rm Aut}^+\mathcal M_i\cong G^+/V_i\cong V_{1-i}\rtimes D^+$ acting primitively and faithfully on its vertices, so by~\cite{Jones} it too is a generalised Paley map, but now one of a chiral pair.

A {\sl generalised Paley map\/} (see~\cite{Jones}) is a Cayley map for the additive group $V$ of a finite field $\mathbb F_q$, where $q=p^d$ for some prime $p$, and where the generating set is a subgroup $S\cong{\rm C}_n$ of the multiplicative group ${\mathbb F}_q^*={\mathbb F}_q\setminus\{0\}$  of ${\mathbb F}_q$ such that
\begin{itemize}
\item $S=-S$ (equivalently $n$ is even or $p$ is odd), and
\item $S$ acts irreducibly on $V$ (equivalently $d=\min\,\{i>0\mid p^i\equiv 1\, {\rm mod}~(n)\}$).
\end{itemize}
Thus the vertex set $V$ is the field ${\mathbb F}_q$. We choose a generator $s$ of the cyclic group $S$, and define the neighbours of each vertex $v\in V$ to be $ v+s, v+s^2,\dots, v+s^n=v+1$ in that cyclic order, following the local orientation around $v$. This defines an orientably regular map ${\mathcal M}(s)$ with ${\rm Aut}^+{\mathcal M}(s)=V\rtimes S\le{\rm AGL}_1({\mathbb F}_q)$. As examples, if $n=q-1$ then $\mathcal M(s)$ is a Biggs embedding of the complete graph $K_q$ (see~\cite{Biggs, JJ85}), while if $q\equiv 1$ mod~$(4)$ and $|S|=(q-1)/2$ it embeds the Paley graph $P_q$. As shown in~\cite{Jones} we have ${\mathcal M}(s)\cong{\mathcal M}(s')$ if and only if $s$ and $s'$are equivalent under the Galois group ${\rm Gal}\,{\mathbb F}_q\cong{\rm C}_d$ of ${\mathbb F}_q$, generated by the Frobenius automorphism $t\mapsto t^p$; thus there are $\phi(n)/d$ such maps for each pair $n$ and $q$, all mutually equivalent under the hole operations $H_j:s\mapsto s^j$ for $j$ coprime to $n$. In particular, the mirror image of $\mathcal M(s)$ is $H_{-1}{\mathcal M}(s)={\mathcal M}(s^{-1})$, so $\mathcal M(s)$ is regular if and only if $s$ and $s^{-1}$ are equivalent under ${\rm Gal}\,{\mathbb F}_q$, that is, $p^e\equiv -1$ mod~$(n)$ for some $e$, necessarily $e=d/2$ with $d$ even.

In case~A we therefore obtain $\phi(n)/d$ regular maps $\mathcal M^+=\mathcal M(s)$, mutually equivalent under hole operations, where $d$ is the multiplicative order of $p$ mod~$(n)$ and $p^{d/2}\equiv -1$ mod~$(n)$. In case~B we have $\dim V_i=\frac{1}{2}\dim V$ for $i=0, 1$, so if $\dim V=d$ then $\dim V_i=e=d/2$ and the maps $V_i$ are a chiral pair of generalised Paley maps over the field ${\mathbb F}_{p^e}$; there are $\phi(n)/e$ such maps, forming $\phi(n)/2e=\phi(n)/d$ chiral pairs, so (as in case~A) we obtain $\phi(n)/d$ maps $\mathcal M^+$, mutually equivalent under hole operations, where now $d=2e$ and $e$ is the multiplicative order of $p$ mod~$(n)$. A more uniform way of describing the two cases is to note that, in either case, $d=2e$ where $e$ is the order of $p$ as an element of the quotient group ${\mathbb Z}_n^*/\langle -1\rangle$, and that the cases are distinguished by whether $p^e\equiv -1$ mod~$(n)$, in case~A, or $p^e\equiv 1$ mod~$(n)$, in case~B.

In case~A one can interpret $\phi(n)/d$ as the number of irreducible factors of the mod~$(p)$ reduction of the cyclotomic polynomial $\Phi_n(t)$, each of degree $d$, giving the possible minimal polynomials on $V$ of a generator of $D^+$. There is a similar interpretation of $\phi(n)/d=\phi(n)/2e$ in case~B, with $e$ replacing $d$; division by $2$ arises from the isomorphisms induced by transposing summands $\mathcal M_i$.

This argument shows that, for a given value of $n$, $p$ must be coprime to $n$, and the structures of $\mathcal M^+$ and $G$, specifically the value of $d$ and whether we are in case~A or B, depend only on the congruence class of $p$ mod~$(n)$, each containing infinitely many primes by Dirichlet's Theorem. To summarise, we have proved the following:

\begin{theo}\label{th:affgps}
If $\mathcal M$ is a regular map of valency $n>2$, and $G:={\rm Aut}\,\mathcal M$ is a permutation group of affine type acting primitively and faithfully on its vertices, then $G=V\rtimes D\le{\rm AGL}_d(p)$ where $V$ is an elementary abelian group of order $q=p^d$ for some prime $p$ coprime to $n$, with $p=2$ if $n$ is odd. Then the vertex set can be identified with $V$, acting regularly by translations, and the subgroup $D:=G_0\cong {\rm D}_n$ of $G$ fixing the vertex $0\in V$ is a subgroup of ${\rm GL}_d(p)$ acting irreducibly on $V$, in one of the following ways, where $D^+$ is the cyclic subgroup of index $2$ in $D$, and $d=2e$ is even:
\begin{enumerate}
\item[(A)] $p$ has multiplicative order $d$ {\rm mod}~$(n)$, with $p^e\equiv -1$ {\rm mod}~$(n)$, and $D^+$ acts irreducibly on $V\cong{\mathbb F}_q$ as a subgroup of ${\mathbb F}_q^*$\,, or
\item[(B)] $p$ has multiplicative order $e$ {\rm mod}~$(n)$, and $V$ decomposes as a $D^+$-module as $V_0\oplus V_1$, with elements of $D\setminus D^+$ transposing the summands $V_i\cong{\mathbb F}_{q'}$ ($q':=\sqrt q=p^e$), and with $D^+$ having mutually inverse irreducible actions on them as a subgroup of ${\mathbb F}_{q'}^*$.
\end{enumerate}

In either case, for a given $n$ and $q$ satisfying the above conditions, there are $\phi(n)/d$ Petrie dual pairs of maps, one map $\mathcal M^+$ orientable and the other map $\mathcal M^-$ non-orientable. In case~A $\mathcal M^+$ is a regular generalised Paley map, embedding the generalised Paley graph $P^{(n)}_q$, and in case B it is the join of a chiral pair of orientably regular generalised Paley maps $\mathcal M_i=\mathcal M^+/V_i$ ($i=0, 1$), each embedding the generalised Paley graph $P^{(n)}_{q'}$.
\end{theo}

 \smallskip

\noindent{\bf Example} For small even $n>2$, the primes $p$ for which $-1$ is a power $p^e$  of $p$ mod~$(n)$, thus giving Case~A, with the corresponding values of $d=2e$, are as follows:
\begin{itemize}
\item for $n=4$, $6$, $8$, $12$, only $p\equiv -1$ mod~$(n)$, with $d=2$;
\item for $n=10$, $p\equiv -1,\pm 3$ mod~$(10)$, with $d=2, 4, 4$;
\item for $n=14$, $p\equiv -1, 3, 5$ mod~$(14)$, with $d=2, 6, 6$.
\end{itemize}
Case~B arises, with $p^e\equiv 1$ mod~$(n)$, in the following cases:
\begin{itemize}
\item for $n=4, 6, 10$, only $p\equiv 1$ mod~$(n)$, with $d=2$;
\item for $n=8$, $p\equiv 1, \pm 3$ mod~$(8)$, with $d=2, 4, 4$;
\item for $n=12$, $p\equiv 1, \pm 5$ mod~$(12)$, with $d=2, 4, 4$;
\item for $n=14$, $p\equiv 1, -3, -5$ mod~$(14)$, with $d=2, 6, 6$.
\end{itemize}
For odd $n$ only the prime $p=2$ corresponds to regular maps; for small $n$ we have:
\begin{itemize}
\item Case~A arises for $n=3, 5, 9, 11, 13$ with $d=2, 4, 6, 10, 12$;
\item Case~B arises for $n=7, 15$ with $d=6, 8$.
\end{itemize}
 
\medskip

\noindent{\bf Remark} The possible dimensions $d$ of $V$ are unbounded. For instance, if we take $n=2m$ for an odd prime $m$ then a generator $c$ for the group ${\mathbb Z}_n^*\cong{\rm C}_{m-1}$ satisfies $c^{(m-1)/2}=-1$ mod~$(n)$, so primes $p\equiv c$ mod~$(n)$ have $d=m-1$ in case~A. If, in addition, we choose $m\equiv -1$ mod~$(4)$ then $c^2$ has odd order $(m-1)/2$, so $-1$ cannot be a power of $c^2$, and hence primes $p\equiv c^2$ mod~$(n)$ have $d=m-1$ in case~B. However, for any fixed $n$ the multiplicative order $d$ or $e$ of $p$ mod~$(n)$, in cases~A or B, must divide $\phi(n)$, so $d$ divides $2\phi(n)$.

\medskip

One can determine the extended type $\{m,n\}_l$, and hence the genus, of each map $\mathcal M^{\pm}$ and $\mathcal M_i$, either by considering the orders of the elements $r_0r_1$ and $r_0r_1r_2$ as in the almost simple case, by using the information about generalised Paley maps in~\cite{Jones}, or from the section on affine groups in~\cite{JLSW}. By~\cite[Lemma~2.4]{Jones}, a generalised Paley map of order $q=p^d$ and valency $n>2$ has type $\{n, n\}_{2p}$ and genus $1+\frac{1}{4}q(n-4)$ unless $n\equiv 2$ mod~$(4)$, in which case it has type $\{n/2,n\}_{2p}$ and genus $1+\frac{1}{4}q(n-6)$. This applies directly to the regular maps $\mathcal M^+$ in case~A. In case~B, replacing $q$ with $q'=\sqrt q$, we see that the chiral pair of maps $\mathcal M_i$ have type $\{n, n\}_{2p}$ and genus $1+\frac{1}{4}q'(n-4)$ unless $n\equiv 2$ mod~$(4)$, in which case they have type $\{n/2,n\}_{2p}$ and genus $1+\frac{1}{4}q'(n-6)$; since the $q'$-sheeted coverings $\mathcal M^+\to\mathcal M_i$ are unbranched, it follows that $\mathcal M^+$ has the same type and genus as in Case~A. It then follows that in either case $\mathcal M^-$ has genus $2+q(np-n-2p)/2p$, and it has type $\{2p,n\}_n$ unless $n\equiv 2$ mod~$(4)$ in which case it has type $\{2p,n\}_{n/2}$.

%%%%%%%%%%%%%

\section{Duality for affine maps}\label{sec:duality}

In the affine case, if $n$ is odd or $n\equiv 0$ mod~$(4)$ then $\mathcal M^+$ has type $\{n,n\}_{2p}$, so there are the possibilities that $\mathcal M^+$ could be self dual or a member of a dual pair. This issue was not considered in~\cite{JLSW} or~\cite{Jones}, so we will deal with it here.

First suppose that $n$ is odd, so that $p=2$. Given standard generators $r_1$ and $r_2$ generating $D$, we need to find a third generators $r_0$ corresponding to the map $\mathcal M^+$. Since $G=V\rtimes D$ we have $r_0=vr$ where $v\in V$ and $r\in D$. Then $r_0^2=1$ if and only if $r$ inverts $v$, $r_0$ commutes with $r_2$ if and only if $r$ and $v$ do, and $G=\langle r_0, r_1, r_2\rangle$ if and only if $v\ne 1$. Since $n$ is odd and $p=2$ these conditions are equivalent to $r=1$ or $r_2$ and $v^r=v\ne 1$. Since $\mathcal M^+$ is orientable we require $r_0\in G^-$, so $r=r_2$ and $v^{r_2}=v\ne 1$.

Now $\mathcal M^+$ is self-dual if and only if there is an automorphism of $G$ fixing $r_1$ and transposing $r_0$ and $r_2$. Conjugation by an element $w\in V$ achieves this if and only if $r_1^w=r_1$ and $r_2^w=r_2v$, or equivalently $w^{r_1}=w$ and $w^{r_2}w=v$. Given $w\in V$ fixed by $r_1$ and not by $r_2$, define $v=w^{r_2}w$, so $v^{r_2}=ww^{r_2}=w^{r_2}w=v$ and we are done. Thus $\mathcal M^+$ is self-dual whenever $n$ is odd.

Now suppose that $n\equiv 0$ mod~$(4)$, so $p$ is odd and coprime to $n$. Then $V$ has mutually coprime order and index, so by the extended Schur--Zassenhaus Theorem it has a single conjugacy class of complements in $G$. Since the faces have valency $n$ their stabilisers are also the vertex stabilisers $D$, so $G$ acts primitively on the faces of $\mathcal M^+$, or equivalently the dual map $D(\mathcal M^+)$ also has automorphism group $G$ acting primitively on its vertices. It follows that either $\mathcal M^+$ is self-dual (as must be the case if $\phi(n)/d=1$), or $\mathcal M^+$ is one of a dual pair of orientable regular maps with this property.

To determine which of these possibilities arises, we can argue as in the case where $n$ is odd, that $r_0=vr$ where $r$ inverts $v$, $v$ and $r$ commute with $r_2$, $v\ne 1$ and $r\in D^-$. Now the elements $r\in D^-$ commuting with $r_2$ are $r_2$ and $r_2z$, where $z:=a^{n/2}$ ($a:=r_1r_2$) generates the centre of $D$. If we take $r=r_2$ we find that $r$ both inverts and fixes $v$, which is impossible since $v$ has order $p>2$. Thus $r=r_2z$, so $r_0=vr_2z$.

The map $\mathcal M^+$ is self-dual if and only if there is an automorphism $\alpha$ of $G$ fixing $r_1$ and transposing $r_0$ and $r_2$. Since $V$ is a characteristic subgroup of $G$, such an automorphism must induce an automorphism $\overline\alpha$ of $G/V\cong{\rm D}_n$ with the same effect on the images $\overline r_i$ of the generators $r_i$ of $G$. Since $\overline\alpha$ fixes $\overline r_1$ and transposes $\overline r_2$ and $\overline r_2\overline z$, it follows that $\overline\alpha$ is the unique automorphism of ${\rm D}_n$ fixing $\overline r_1$ and sending $\overline a$ to $\overline a\overline z$.

We need to determine whether or not this automorphism $\overline\alpha$ of ${\rm D}_n$ lifts to an automorphism $\alpha$ of $G$ with the required properties, that is, fixing $r_1$ and transposing $r_2$ and $r_0=vr_2z$.
This happens if and only if the representation of $D$ on $V$ is equivalent to its composition with $\overline\alpha$ (applied to $D$); since $z$ is represented on $V$ (in case~A) or each $V_i$ (in case~B) by the matrix $-I$, this is equivalent to the eigenvalues $\lambda$ of $a$ being Galois conjugate to $-\lambda$ or $-1/\lambda$. Since $\lambda$ has order $n$ we see that if $n\equiv 0$ mod~$(4)$ then $\mathcal M^+$ is self-dual if and only if $p^i\equiv \frac{n}{2}\pm1$ mod~$(n)$ for some $i$.

%%%%%%%%%%%%

\section{Affine examples for small $n$}\label{sec:affexs}

We first consider odd values of $n$, so $p=2$ and the maps $\mathcal M^+$ and $\mathcal M^-$ have genus $1+2^{d-2}(n-4)$ and $2+2^{d-2}(n-4)$.

\smallskip

{\bf Example} If $n=3$ we have case~A with $d=2$ and $G\cong{\mathbb F}_4\rtimes{\rm D}_3\cong{\rm S}_4$. The number of Petrie dual pairs $\mathcal M^{\pm}$ is $\phi(3)/d=2/2=1$. The map $\mathcal M^+$ is the tetrahedron ${\mathcal M}^+=\{3,3\}$, of genus $0$ and type $\{3,3\}_4$, while ${\mathcal M}^-$, of non-orientable genus $1$ and type $\{4,3\}_3$, is the antipodal quotient of the cube. For each map, $G$ acts $4$-transitively (and thus primitively) on the vertices.

\medskip

{\bf Example} If $n=5$ we have case~A with $d=4$ and $G\cong{\mathbb F}_{16}\rtimes{\rm D}_5$. The unique Petrie dual pair consists of ${\mathcal M}^+=$ R5.9 (self-dual) of type $\{5,5\}_4$ and ${\mathcal M}^-=$ N6.3 of  type $\{4,5\}_5$.

\medskip

{\bf Example} If $n=7$ we have case~B with $d=6$ and $G\cong({\mathbb F}_8\oplus{\mathbb F}_8)\rtimes{\rm D}_7$. There is a unique Petrie dual pair, consisting of R49.57 (self-dual) of type $\{7,7\}_4$ and N50.5 of type $\{4,7\}_7$. The first of these is the join of the chiral pair C7.2 of Edmonds embeddings ${\mathcal M}_i\;(i=1, 2)$ of the complete graph $K_8$, each with automorphism group $G^+/V_i\cong{\rm AGL}_1(8)\cong{\mathbb F}_8\rtimes{\rm C}_7$. (See~\cite{JJ85} for orientably regular embeddings of complete graphs.)
\medskip

{\bf Example} If $n=9$ we have case~A with $d=6$ and $G\cong{\mathbb F}_{64}\rtimes{\rm D}_9$. The unique Petrie dual pair consists of ${\mathcal M}^+=$ R81.125 (self-dual) of type $\{9,9\}_4$ and ${\mathcal M}^-=$ N82.1 of type $\{4,9\}_9$.

\medskip

{\bf Example} If $n=11$ we have case~A with $d=10$ and $G\cong{\mathbb F}_{1024}\rtimes{\rm D}_{11}$. The unique Petrie dual pair are ${\mathcal M}^+$ (self-dual) of genus $1793$ and type $\{11,11\}_4$ and ${\mathcal M}^-$ of genus $1794$ and type $\{4,11\}_{11}$. For odd $n\ge 11$ the genera are too large for the maps to appear in~\cite{Conder}.

\medskip

{\bf Example} If $n=13$ we have case~A with $d=12$ and $G\cong{\mathbb F}_{4096}\rtimes{\rm D}_{13}$. There is a unique Petrie dual pair, namely ${\mathcal M}^+$ (self-dual) of genus $9217$ and type $\{13,13\}_4$, and ${\mathcal M}^-$ of genus $9218$ and type $\{4,13\}_{13}$.

\medskip

{\bf Example} If $n=15$ we have case~B with $d=8$ and $G\cong({\mathbb F}_{16}\oplus{\mathbb F}_{16})\rtimes{\rm D}_{15}$. The unique Petrie dual pair are $\mathcal M^+$ of genus $705$ and type $\{15,15\}_4$ (self-dual), and $\mathcal M^-$ of genus $706$ and type $\{4,15\}_{15}$. The first of these is the join of the chiral pair C45.2 of orientably regular embeddings of $K_{16}$, each with automorphism group ${\rm AGL}_1(16)\cong{\mathbb F}_{16}\rtimes{\rm C}_{15}$. 

\medskip

In each of the above examples the number $\phi(n)/d$ of Petrie dual pairs ${\mathcal M}^{\pm}$ is equal to $1$. This will always happen if $2$ is a primitive root mod~$(n)$, that is, a generator of the group ${\mathbb Z}_n^*$, a condition which implies that $n=r^i$ where $r$ is prime and $2$ is a primitive root mod~$(r)$. Indeed, if $r$ is prime and $2$ is a primitive root mod $r^2$, as is the case for $r=3, 5, 11$ and $13$, for instance, then $2$ is a primitive root mod~$(r^i)$ for all $i$ (see~\cite[Lemma~1.10(i)]{Nar}), giving infinitely many $n$ with $\phi(n)/d=1$. Artin's Primitive Roots Conjecture asserts that $2$ (along with most other integers) is a primitive root for a set of primes of positive density. However, it is still an open problem whether there are infinitely many such primes.

In contrast to the above examples, there are two Petrie dual pairs for $n=17$, and the following examples show that there can be many more than that.

\medskip

{\bf Example} A simple argument using the multiplicative property of Euler's function shows that $\phi(n)\ge \sqrt n/2$ for all $n\in{\mathbb N}$. If $n$ is a Mersenne number $2^e-1$ then $2$ has order $e$ mod~$(n)$, giving a lower bound $\phi(n)/d=\phi(n)/2e\ge \sqrt{2^e-1}/4e$, so the number of Petrie dual pairs is unbounded above. In fact, if $e>2$ then $-1$ is not a power of $2$ mod~$(n)$, so we have case~B with $G\cong({\mathbb F}_{2^e}\oplus{\mathbb F}_{2^e})\rtimes{\rm D}_{n}$. Each map $\mathcal M^+$ is the join of a chiral pair of orientably regular embeddings $\mathcal M_i$ of $K_{2^e}$, each of genus  $1+2^{e-2}(n-4)$ with automorphism group ${\rm AGL}_1(2^e)\cong{\mathbb F}_{2^e}\rtimes{\rm C}_n$. For instance, if $e=5$ and $n=31$ there are three Petrie dual pairs $\mathcal M^{\pm}$ of genus $6913$ and $6914$; the three chiral pairs $\mathcal M_i$ are the orientably regular embeddings  C217.45, C217.46 and C217.47 of $K_{32}$.

Similarly we have an infinite family of valencies $n = 2^e + 1$ with an unbounded number of vertex-primitive regular embeddings of generalised Paley graphs.
Likewise, for any odd prime $p$ the sequences $p^e\pm 1$ provide valencies with an unbounded number of Petrie dual maps in both Case~A and Case~B.

\medskip

{\bf Example} If $n=p_1\ldots p_k$ where the $p_i$ are distinct primes for which $2$ is not a primitive root (for instance if $p_i\equiv \pm 1$ mod~$(8)$, since $2$ is then a square mod~$(p)$ by Euler's Criterion), then $\phi(n)/d\ge 2^k$; there are infinitely many such primes $p_i$ by Dirichlet's Theorem, so the number of Petrie dual pairs is unbounded above.

\medskip

We now consider examples where $n$ is even. In these cases $p$ can be any prime not dividing $n$.

\medskip

{\bf Example} Let $n=4$. If $p\equiv 1$ mod~$(4)$ we have case~B with $d=2$ and $G\cong({\mathbb F}_p\oplus{\mathbb F}_p)\rtimes{\rm D}_4$, whereas if $p\equiv -1$ mod~$(4)$ we have case~A with $d=2$ and $G\cong{\mathbb F}_{p^2}\rtimes{\rm D}_4$. In either case the number of Petrie dual pairs is $\phi(4)/d=1$. They are the self-dual torus map ${\mathcal M}^+=\{4,4\}_{p,0}$ and its Petrie dual ${\mathcal M}^-$, a non-orientable map of genus $(p-1)^2+1$ and type $\{2p,4\}_4$. For instance, if $p=3, 5, 7$ or $11$ then $\mathcal M^-$ is N5.2, N17.2, N37.1 or N101.3. If $p\equiv 1$ mod~$(4)$ then $p=a^2+b^2$ for integers $a, b>0$, and $\mathcal M^+$ is the join of the chiral pair $\mathcal M_i=\{4,4\}_{(a,b)}$ and $\{4,4\}_{(b,a)}$ of orientably regular torus maps; if $p=5$ or $13$ these are embeddings of $K_5$ or the generalised Paley graph $P_{13}^{(4)}$ of order $13$ and valency $4$ (see~\cite{Jones} for orientably regular embeddings of Paley and generalised Paley graphs).

\medskip

{\bf Example} When $n=6$ the situation is similar. For primes $p\equiv \pm 1$ mod~$(6)$ we have case~B or A respectively. The two maps are the torus map ${\mathcal M}^+=\{3,6\}_{p,0}$ and its Petrie dual ${\mathcal M}^-$, a non-orientable map of type $\{2p,6\}_3$ and genus $2p^2-3p+2=2(p-1)^2+p$. For instance, if $p=5, 7, 11$ or $13$ this is N37.5, N79.2, N211.2 or N301.2. If $p\equiv 1$ mod~$(6)$ then $p=a^2+ab+b^2$ for integers $a, b>0$, and $\mathcal M^+$ is the join of the chiral pair $\mathcal M_i=\{3,6\}_{(a,b)}$ and $\{3,6\}_{(b,a)}$ of orientably regular torus maps; if $p=7$ or $13$ these are embeddings of $K_7$ or the Paley graph $P_{13}$.

\medskip

{\bf Example} Let $n=8$. If $p\equiv \pm 1$ mod~$(8)$ we have case~B or A respectively, with $d=2$, $G\cong({\mathbb F}_p\oplus{\mathbb F}_p)\rtimes{\rm D}_8$ or ${\mathbb F}_{p^2}\rtimes{\rm D}_8$, and $\phi(8)/d=4/2=2$ Petrie dual pairs. Each map $\mathcal M^+$ has genus $p^2+1$ and type $\{8,8\}_{2p}$, while $\mathcal M^-$ has genus $3p^2-4p+2$ and type $\{2p,8\}_8$. If $p=7$, for instance, there are a dual pair R50.7 of maps ${\mathcal M}^+$ of type $\{8,8\}_{14}$, and their Petrie duals N121.3 and N121.4 of type $\{14,8\}_8$. If $p=17$ there are a dual pair R290.6 of maps $\mathcal M^+$ of type $\{8,8\}_{34}$, together with their Petrie duals $\mathcal M^-$ of genus $801$ and type $\{34,8\}_8$, the latter having too large genus to appear in~\cite{Conder}. In this case the four quotient maps $\mathcal M_i$ are a dual pair C18.1 of chiral pairs of orientably regular maps of type $\{8,8\}_{34}$, each of them embedding the Paley graph $P_{17}$.

If $p\equiv \pm 3$ mod~$(8)$ we have case~B with $d=4$, $G\cong({\mathbb F}_{p^2}\oplus{\mathbb F}_{p^2})\rtimes{\rm D}_8$ and just one Petrie dual pair, $\mathcal M^+$ of genus $p^4+1$ and type $\{8,8\}_{2p}$ and $\mathcal M^-$ of genus $3p^4-4p^3+2$ and type $\{2p,8\}_8$. For instance, if $p=3$ then $\mathcal M^+$ is the self-dual map R82.50 of type $\{8,8\}_6$, and $\mathcal M^-$ is its Petrie dual N137.3 of type $\{6,8\}_8$; $\mathcal M^+$ is the join of the chiral pair C10.3 of self-dual orientably regular embeddings of $K_9$. If $p=5$ then $\mathcal M^+$ has genus $626$ and type $\{8,8\}_{10}$, while $\mathcal M^-$ has genus $1377$ and type $\{10,8\}_8$; $\mathcal M^+$ is the join of the chiral pair C26.1 of mutually dual orientably regular embeddings of the generalised Paley graph $P_{25}^{(8)}$.

\medskip

{\bf Example} The situation is similar for $n=10$. If $p\equiv \pm 1$ mod~$(10)$ we have case~B or A respectively, with $d=2$, $G\cong({\mathbb F}_p\oplus{\mathbb F}_p)\rtimes{\rm D}_{10}$ or ${\mathbb F}_{p^2}\rtimes{\rm D}_{10}$, and $\phi(10)/d=2$ Petrie dual pairs. Each map $\mathcal M^+$ has genus $p^2+1$ and type $\{5,10\}_{2p}$, while $\mathcal M^-$ has genus $4p^2-5p+2$ and type $\{2p,10\}_5$. For instance, if $p=11$ then $\mathcal M^+$ is R122.3 or R122.4 of type $\{5,10\}_{22}$, while its Petrie dual $\mathcal M^-$ is respectively N431.2 or N431.1 of type $\{22,10\}_5$; the maps $\mathcal M^+=$ R122.3 and R122.4 are respectively the joins of the two chiral pairs C12.2 and C12.1 of orientably regular embeddings of $K_{11}$. If $p=19$ the maps $\mathcal M^+$ have genus $362$ and type $\{5,10\}_{38}$, while their Petrie duals $\mathcal M^-$ have genus $1361$ and type $\{38,10\}_5$.

If $p\equiv \pm 3$ mod~$(10)$ we have case~A with $d=4$, $G\cong {\mathbb F}_{p^4}\rtimes{\rm D}_{10}$, and one Petrie dual pair.
Each map $\mathcal M^+$ has genus $p^4+1$ and type $\{5,10\}_{2p}$, while $\mathcal M^-$ has genus $4p^2-5p+2$ and type $\{2p,10\}_5$. For instance, if $p=3$ the maps are R82.19 of type $\{5,10\}_6$ and N191.7 of type $\{6,10\}_5$. If $p=7$ then $\mathcal M^+$ has genus $2402$ and type $\{5,10\}_{14}$, while $\mathcal M^-$ has genus $7891$ and type $\{14,10\}_5$.

\medskip

{\bf Example} Let $n=12$. If $p\equiv\pm 1$ mod~$(12)$ we have case~B or A with $d=2$, $G\cong ({\mathbb F}_p\oplus{\mathbb F}_p)\rtimes{\rm D}_{12}$ or ${\mathbb F}_{p^2}\rtimes{\rm D}_{12}$, and $\phi(12)/d=2$ Petrie dual pairs. The maps $\mathcal M^+$ have genus $2p^2+1$ and type $\{12,12\}_{2p}$, while their Petrie duals $\mathcal M^-$ have genus $5p^2-6p+2$ and type $\{2p,12\}_{12}$. For instance, if $p=11$ there are four maps, consisting of a dual pair R243.16 of type $\{12,12\}_{22}$ and their Petrie duals N541.6 and N541.7 of type $\{22,12\}_{12}$. (In~\cite{Conder} one also finds a Petrie dual pair R243.17 and N541.8 with the same extended types and genera as these four maps.  However their automorphism group, also a semidirect product ${\mathbb F}_{11^2}\rtimes{\rm D}_{12}$, is not isomorphic to that of the other four: it has a primitive but non-faithful action on the vertices, with a kernel of order $3$, showing that R243.17 is a triple cover of the regular torus map $\{4,4\}_{11,0}$, branched over the vertices and face-centres. See Section~\ref{sec:non-faithful} for this example and for more on non-faithful actions.)
If $p=13$ each map $\mathcal M^+$ has genus $339$ and type $\{12,12\}_{26}$, while their Petrie duals $\mathcal M^-$ have genus $769$ and type $\{26,12\}_{12}$; the four quotient maps $\mathcal M_i$ are the two dual chiral pairs C27.7 of orientably regular embeddings of $K_{13}$.

If $p\equiv\pm 5$ mod~$(12)$ we have case~B with $d=4$, $G\cong ({\mathbb F}_{p^2}\oplus{\mathbb F}_{p^2})\rtimes{\rm D}_{12}$, and one Petrie dual pair, namely $\mathcal M^+$ of genus $2p^4+1$ and type $\{12,12\}_{2p}$ and $\mathcal M^-$ of genus $5p^4-6p^3+2$ and type $\{2p,12\}_{12}$. If $p=5$ then $\mathcal M ^+$ has genus $1251$ and type $\{12,12\}_{10}$, while $\mathcal M^-$ has genus $2377$ and type $\{10,12\}_{12}$. The two quotient maps $\mathcal M_i$ are the chiral pair C101.31 of orientably regular self-dual embeddings of type $\{12,12\}_{10}$ of the Paley graph $P_{25}$, so $\mathcal M^+$ is also self-dual. 

%%%%%%%%%%%%%%%%%%%

%\iffalse % Delete thissection for the SIGMAP paper but retain for the arXiv?

\section{Dihedral quotient maps}\label{sec:dihedral}

%{\color{blue}[This section was omitted from the ADAM version. Include it here?]}

%\medskip

In the affine case $G$ has a normal subgroup $V$ with $\overline G:=G/V\cong {\rm D}_n$, so each $\mathcal M^{\pm}$ is a regular covering, by $V$, of a regular map $\overline{{\mathcal M}}:={\mathcal M}/V$ with automorphism group $\overline G$. Now $\mathcal M^+$ has type $\{m, n\}$ with $m=n/2$ or $n$ as $n\equiv 2$ mod~$(4)$ or not; since $m$ and $n$ are coprime to the exponent $p$ of $V$, it follows that $\overline{\mathcal M^+}$ has the same type $\{m,n\}$ as $\mathcal M^+$, and has two $m$-valent faces or one, respectively. If $n$ is even the covering $\mathcal M\to\overline{\mathcal M^+}$ is unbranched, so $\overline{\mathcal M^+}$ has genus $\overline{g}=1+(g-1)/q=n/4$ or $(n-2)/4$ as $\mathcal M^+$ has genus $g=1+q(n-4)/4$ or $1+q(n-6)/4$ for $n\equiv 0$ or $2$ mod~$(4)$; in fact, $\overline{\mathcal M^+}$ is the dual of the familiar picture of an $n$-gon with opposite sides identified to form an orientable surface (see Figure~\ref{n=4,6} for $n=4$ and $6$).

\begin{figure}[h!]
\label{fig:cube}
\begin{center}
\begin{tikzpicture}[scale=0.15, inner sep=0.7mm]

\node (a) at (10,10)  {};
\node (b) at (-10,10) {};
\node (c) at (-10,-10) {};
\node (d) at (10,-10) {};
\node (e) at (0,0) [shape=circle, fill=black] {};

\draw [thick, dashed] (a) to (b) to (c) to (d) to (a);
\draw [thick] (-10,0) to (10,0);
\draw [thick] (0,-10) to (0,10);

%%%%%%%%%%

\node (O) at (40,0) [shape=circle, fill=black] {};
\node (A) at (53,0)  {};
\node (B) at (46,10)  {};
\node (C) at (34,10)  {};
\node (D) at (27,0)  {};
\node (E) at (34,-10)  {};
\node (F) at (46,-10)  {};
\draw [thick, dashed] (A) to (B) to (C) to (D) to (E) to (F) to (A);
\draw [thick] (40,-10) to (40,10);
\draw [thick] (30.5,-5) to (49.5,5);
\draw [thick] (30.5,5) to (49.5,-5);

\end{tikzpicture}

\end{center}
\caption{Quotient maps $\overline{\mathcal M^+}$ for $n=4$ and $n=6$ (opposite sides identified)}
\label{n=4,6}
\end{figure}
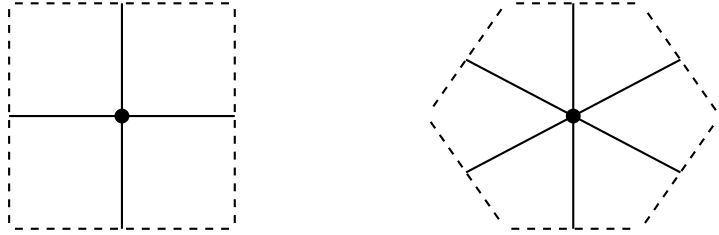

 If $n$ is odd then $r_0r_2\in V$ and the covering $\mathcal M\to\overline{\mathcal M^+}$ is branched over the edges, so that $\overline{\mathcal M^+}$ is a map on the sphere with one vertex, one face, and $n$ half-edges. There is also branching of order $p$ over the Petrie polygons, so that whereas $\mathcal M^+$ has extended type $\{m,n\}_{2p}$, $\overline{\mathcal M^+}$ has extended type $\{m,n\}_2$. Taking Petrie duals, we see that $\overline{\mathcal M^-}$ has extended type $\{2,n\}_m$; when $n$ is even this map is the antipodal quotient of the $n$-valent beachball map $\{2,n\}$ on the sphere, and when $n$ is odd it is a map on the closed disc, with one vertex at the centre and $n$ half-edges from there to the boundary.
 %{\color{blue}[Add an analogue of Figure~\ref{n=4,6} to illustrate these maps?]}

Since $V$ is an elementary abelian $p$-group, one can obtain the $\phi(n)/d$ original maps $\mathcal M^+$ from $\overline{\mathcal  M^+}$ by taking its mod~$(p)$ homology cover (the `Macbeath trick') and decomposing the covering group, the homology group $H_1(\overline{\mathcal M^+};{\mathbb Z}_p)\cong{\mathbb Z}_p^{2\overline g}$, as a direct sum of irreducible $D$-submodules. Each map $\mathcal M^+$ is then obtained by factoring out the sum of all but one of these submodules.

%\fi

%%%%%%%%%%%%%%%%

\section{Non-faithful actions}\label{sec:non-faithful}

Regular maps with automorphism groups acting non-faithfully on vertices, edges or faces have been dealt with thoroughly by Li and \v Sir\'a\v n in~\cite{LS}, so here we will just outline some possibilities.

Suppose that $G={\rm Aut}\,{\mathcal M}$ for some regular map $\mathcal M$, acting primitively but not faithfully on the vertices. As before, the vertex stabilisers $D$ are dihedral maximal subgroups of $G$, but now their intersection $K$, the kernel of this action, is nontrivial. Since $K$ is normal in $G$ it is normal in $D$, so either $K\le D^+$, or $n$ is even and $K\cong{\rm D}_{n/2}$. In the latter case either $r_1$ or $r_2$ is an element of $K$, so $\mathcal M$ has non-empty boundary, a case we are excluding. Hence we may assume that $K\le D^+$, so $\mathcal M$ is a regular cyclic covering of the regular map $\overline{\mathcal M}={\mathcal M}/K$. The latter has automorphism group ${\overline G}=G/K$ acting primitively and faithfully on its vertices, so $\overline{\mathcal M}$ and $\overline G$ are as described earlier, either almost simple or affine. The action by conjugation of $G$ on $K$ gives a homomorphism $\overline G\to{\rm Aut}\,K$.
If $|K|>2$ this homomorphism is non-trivial since $r_1$ and $r_2$ invert $K$; now ${\rm Aut}\,K$ is abelian since $K$ is cyclic, so $\overline G$ cannot be perfect.

If $\overline G$ is almost simple this excludes all cases apart from $\overline G={\rm PGL}_2(q)$ for odd $q>5$ (cases~(2) and (5) of Theorem~\ref{th:LiThm}). This group can also be excluded since we have seen that, whereas its standard generators $\overline{r_1}$ and $\overline{r_2}$ both invert $\overline{D^+}$, one of them must lie in the commutator subgroup $\overline{G^+}={\rm PSL}_2(q)$ of $\overline G$, which acts trivially on $K$. This leaves the case $|K|=2$, where $n$ is even and $K=\langle z\rangle$ with $z=(r_1r_2)^{n/2}$, the central involution in $D$.

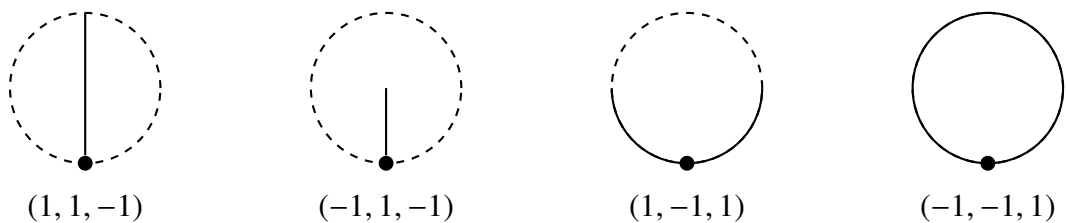
\begin{figure}[h!]
\label{fig:cube}
\begin{center}
\begin{tikzpicture}[scale=0.2, inner sep=0.7mm]

\draw [thick, dashed] (-30,0) arc (0:360:5);
\draw [thick, dashed] (-10,0) arc (0:360:5);
\draw [thick, dashed] (10,0) arc (0:360:5);
\draw [thick, dashed] (30,0) arc (0:360:5);

\node (a) at (-35,-5) [shape=circle, fill=black] {};
\node (b) at (-15,-5) [shape=circle, fill=black] {};
\node (c) at (5,-5) [shape=circle, fill=black] {};
\node (d) at (25,-5) [shape=circle, fill=black] {};

\draw [thick] (a) to (-35,5);
\draw [thick] (b) to (-15,0);
\draw [thick] (10,0) arc (0:-180:5);
\draw [thick] (30,0) arc (0:360:5);

\node at (-35,-8) {$(1,1,-1)$};
\node at (-15,-8) {$(-1,1,-1)$};
\node at (5,-8) {$(1,-1,1)$};
\node at (25,-8) {$(-1,-1,1)$};

\end{tikzpicture}

\end{center}
\caption{The four $2$-flag maps $\mathcal T$ with $r_1r_2\ne 1$}
\label{4mapsT}
\end{figure}

One can construct such double coverings by taking $\mathcal M$ to be the join of $\overline{\mathcal M}$ with a $2$-flag map $\mathcal T$ not covered by $\overline{\mathcal M}$, so that $G=\overline{G}\times{\rm C}_2$. There are seven $2$-flag maps $\mathcal T$, corresponding to the seven generating triples $(r_i)$ for ${\rm Aut}\,{\mathcal T}\cong{\rm C}_2$, and using each has the effect of replacing one, two or all three of the parameters $\overline l, \overline m$ and $\overline n$ in the extended type $\{\overline m, \overline n\}_{\overline l}$ of $\overline{\mathcal M}$ with $l:={\rm lcm}\{\overline l,2\}$, and similarly for $m$ and $n$. In particular, four of the seven maps $\mathcal T$, those with $r_1r_2\ne1$, replace $\overline n$ with $n$, so that if $\overline n$ is odd we obtain a double covering ${\mathcal M}\to\overline{\mathcal M}$ which is bijective on the vertices, doubling their valencies. These maps $\mathcal T$, all on the closed disc, are shown in Figure~\ref{4mapsT}, together with the corresponding generating triples $(r_i)$ for ${\rm Aut}\,{\mathcal T}$, identified with $\{\pm 1\}$.

\medskip

{\bf Example} If $\overline{\mathcal M}$ is the antipodal quotient of the icosahedron, of type $\{3,5\}_5$ with primitive and faithful automorphism group $\overline G\cong{\rm SL}_2(4)\cong{\rm A}_5\cong{\rm PSL}_2(5)$ (see case~(6) in Section~\ref{sec:countingAS}), then the four double covers $\mathcal M$ arising in this way are N6.1 and N6.2, of types $\{3,10\}_{10}$ and $\{3,10\}_5$, branched over the six vertices of $\overline{\mathcal M}$, and N16.5 and N16.6, of types $\{6,10\}_5$ and $\{6,10\}_{10}$, branched over the six vertices and ten face-centres. Each of these four maps has automorphism group $G\cong\overline{G}\times{\rm C}_2$ acting primitively but non-faithfully on its six vertices. (The other three double covers of $\overline{\mathcal M}$ arising as joins are N10.4 and N10.5 of types $\{6,5\}_{10}$ and $\{6,5\}_5$, branched over the face-centres, and the icosahedron itself, which is an unbranched cover; none of these is vertex-primitive.) By taking Petrie duals of $\overline{\mathcal M}$ and $\mathcal M$ we obtain four further examples. The only other double cover of $G$ is the binary icosahedral group $\hat{\rm A}_5\cong{\rm SL}_2(5)$, but since this contains only one involution (generating the centre), it does not yield any regular maps.

\medskip

In contrast, affine groups give many examples, which can be created as follows. Start with an affine group $\overline G=V\rtimes_{\phi} \overline D$ acting  primitively and faithfully on the vertex set $V$ of a regular map $\overline{\mathcal M}$, with a specific action $\phi:\overline D\to {\rm Aut}\,V={\rm GL}_d(p)$, then take a dihedral group $D$ with a normal subgroup $K$ of odd order such that $D/K\cong\overline D$, and consider the covering group $G=V\rtimes_{\psi} D$ where $\psi$ is the composition of $\phi$ with the natural epimorphism $D\to\overline D$.

\medskip

{\bf Example }Let $\overline{\mathcal M}$ be the torus map $\{4,4\}_{(p,0)}$ for some odd prime $p$, with primitive and faithful automorphism group $\overline G=V\rtimes \overline D$ where $V\cong{\rm C}_p\times{\rm C}_p$ and $\overline D\cong{\rm D}_4$. Now take $D\cong{\rm D}_{4k}$ for some odd $k>1$, and let $G=V\rtimes D$ where $D$ acts on $V$ as $\overline D$ via the natural epimorphism $D\to\overline D$ with kernel $K\cong {\rm C}_k$. Then $G$ is the automorphism group of a vertex-primitive orientable regular map $\mathcal M$ which is a $k$-sheeted covering of $\overline{\mathcal M}$, branched over the vertices and face-centres; it is self-dual, of type $\{4k,4k\}_{2p}$ and genus $(k-1)p^2+1$. The action of $G$ on the vertices (and on the faces) has kernel $K$. For instance, if $k=3$ then for $p=5,7$ and $11$ we obtain the maps R51.24, R99.23 and R243.17.

However, when $p=k$ there are $\phi(p)=p-1$ such vertex-primitive but non-faithful covers of $\overline{\mathcal M}$, two of them self-dual and $(p-3)/2$ dual pairs. They are determined by the $(p-1)^2$ choices for a pair of monodromy permutations $\pi_v$ and $\pi_f$, each generating $K\cong{\rm C}_k\cong{\mathbb Z}_p$, at the vertices and at the faces of $\overline{\mathcal M}$; two covers are isomorphic if and only if they have the same ratio $\rho=\pi_v/\pi_f\in{\mathbb Z}_p$, and duality corresponds to inverting $\rho$, so that the self-dual covers are those with $\rho=\pm 1$. For instance, if $p=k=5$ we have self-dual covers R101.41 and R101.43 and a dual pair R101.42;  for $p=k=7$ we have self-dual covers R295.55 and R295.56 and dual pairs R295.57 and R295.58.

%%%%%%%%%%%%%%%%%

\section{Hypermaps}\label{sec:hypermaps}

Almost everything is the same for hypermaps as for maps, with the fact that the vertex stabilisers $D=\langle r_1, r_2\rangle$ are dihedral maximal subgroups again implying that $G$ is either almost simple or affine. However we lose the relation $r_0r_2=r_2r_0$, so instead of looking for involutions $r_0\in C_G(D)\setminus D$ we look for them in $G\setminus D$, giving many more hypermaps than maps. A paper addressing this in more detail is in preparation; here we just give a few simple examples.

\medskip

\noindent {\bf Example} In the almost simple case (6), with $q=4$, $n=5$ and $G={\rm SL}_2(4)\cong{\rm A}_5$ there are  two maps and six proper hypermaps (not maps, up to duality or triality). The latter appear in~\cite{Conder} as NPH6.1 (one hypermap, of type $(5,3,3)$), NPH10.2, NPH10.3 (two each, of types $(5,3,5)$ and $(5,5,3)$) and NPH14.1 (one, of type $(5,5,5)$).

\medskip

Also, unlike in the case of maps, there are affine examples with $n$ and $p$ both odd. This is because, in the case of maps, the argument that $p=2$ if $n$ is odd requires that $(r_0r_2)^2=1$.

\medskip

\noindent {\bf Example} There is an orientable regular hypermap $\mathcal H$ of type $(3,p,3)$ and genus $(p-1)(p-2)/2$ for each prime $p>3$, with $G\cong({\mathbb F}_p\oplus{\mathbb F}_p)\rtimes{\rm D}_3$ or ${\mathbb F}_{p^2}\rtimes{\rm D}_3$ as $p\equiv 1$ or $-1$ mod~$(3)$. If $p\equiv 1$ mod~$(3)$ then $\mathcal H$ is the minimal regular cover of a chiral pair of orientably regular hypermaps $\mathcal H_i$ ($i=1, 2$) of type $(3,p,3)$ and genus~$(p-1)/2$; for instance, when $p=7$ these are the hypermaps CH3.1 of type $(3,7,3)$ on Klein's quartic curve, with regular cover RPH15.1. If $p=5$ then $\mathcal H$ is RPH6.1.

%%%%%%%%%%%%%%%%%%

\section{A useful congruence}\label{sec:congruence}

The affine maps and automorphism groups we have discussed are divided structurally into two cases, A and B, according to whether or not $p^e\equiv -1$ mod~$(n)$ for some $e$ (see Section~\ref{sec:affine} and in particular Theorem~\ref{th:affgps}). Here we characterise the pairs $n, p$ satisfying this condition.

Let $n_r=r^{e_r}$ ($r$ prime, $e_r\ge 1$) be the prime powers appearing in the factorisation of $n$, and let $d_r$ be the multiplicative order of $p$ mod~$(n_r)$. Then we have the following:

\begin{prop}\label{prop:congruence}
We have $p^e\equiv-1$ {\rm mod}~$(n)$ for some $e$ if and only if the orders $d_r$ for primes $r>2$ dividing $n$ all have the same $2$-part $2^k\ge2$, with $p\equiv-1$ {\rm mod}~$(n_2)$ if $n_2\ge 2$, and $k=1$ if $n_2\ge 4$.
%In these circumstances the congruence $p^i\equiv-1$ {\rm mod}~$(n)$ is satisfied by $i={\rm lcm}\,\{d_r/2\mid r>2\}$ if some prime $r>2$ divides $n$, and otherwise by $i=1$.
\end{prop}

The proof is straightforward, using the fact that the group of units mod~$(n_r)$ is cyclic unless $r=2$ and $n_2\ge 8$, in which case it has the form $\langle 5\rangle\times\langle-1\rangle\cong{\rm C}_{n_2/4}\times{\rm C}_2$.
The primality of $p$ is irrelevant here: the conditions depend only on the congruence class of $p$ mod~$(n)$. Of course, by Dirichlet's Theorem every congruence class in ${\mathbb Z}_n^*$ contains infinitely many primes.

\section{Acknowledgments}

The second author acknowledges support from the APVV Research
Grant APVV-19-0308 and from the VEGA Research Grants 1/0423/20 and 1/0727/22.

%%%%%%%%%%%%%%%%
%%%%%%%%%%%%%%%%

% end{document}

%%%%%%%%%%%%%%%%%%%%%%%%%%%%%
%%%%%%%%%%%%%%%%%%%%%%%%%%%%%

%\iffalse % Omit Appendices A and B for the ADAM version.

%\newpage

\section{Appendix A: data for maps with almost simple groups}

%{\color{blue}[Appendices A and B were omitted from the ADAM version.
%They should be included in the arXiv version, possibly with labels and captions added to the tables so that we can add references to them in the main text.
%There is also room for more tables, if required.]}

%\medskip

Here, for each of cases $i=1,\ldots, 7$ where the automorphism group $G$ is almost simple, we give the number $\nu_i(q)$ of vertex-primitive maps $\mathcal M$, together with their minimum and maximum genus, for each relevant prime power $q<50$ (see Section~\ref{sec:countingAS}).

$$\begin{array}{|c|c|c|c|c|c|}\hline
\mbox{Case}& q & n & \mbox{Min.~genus} & \mbox{Max.~genus} & \nu_i(q) \\ \hline
(1)  &  13  &  6  &  106  &  142  &  2  \\ \hline
(1)  &  17  &  8  &  325  &  389  &  6  \\ \hline
(1)  &  19  &  9  &  97  &  496  &  15  \\ \hline
(1)  &  23  &  11  &  232  &  991  &  30  \\ \hline
(1)  &  25  &  12  &  847  &  1327  &  5  \\ \hline
(1)  &  27  &  13  &  1379  &  1730  &  14  \\ \hline
(1)  &  29  &  14  &  582  &  2402  &  18  \\ \hline
(1)  &  31  &  15  &  746  &  2761  &  32  \\ \hline
(1)  &  37  &  18  &  1408  &  5284  &  24  \\ \hline
(1)  &  41  &  20  &  2011  &  7331  &  36  \\ \hline
(1)  &  43  &  21  &  2367  &  8086  &  66  \\ \hline
(1)  &  47  &  23  &  3198  &  10765  &  132  \\ \hline
(1)  &  49  &  24  &  3677  &  12301  &  22  \\ \hline
\end{array}$$

%%%%%%%%%%%%%%

$$\begin{array}{|c|c|c|c|c|c|}\hline
\mbox{Case}& q & n & \mbox{Min.~genus} & \mbox{Max.~genus} & \nu_i(q) \\ \hline
(2)  &  7  &  6  &  16  &  37  &  5  \\ \hline
(2)  &  9  &  8  &  17  &  101  &  7  \\ \hline
(2)  &  11  &  10  &  46  &  211  &  18  \\ \hline
(2)  &  13  &  12  &  93  &  379  &  22  \\ \hline
(2)  &  17  &  16  &  257  &  937  &  60  \\ \hline
(2)  &  19  &  18  &  382  &  1351  &  51  \\ \hline
(2)  &  23  &  22  &  738  &  2509  &  105  \\ \hline
(2)  &  25  &  24  &  977  &  3277  &  46  \\ \hline
(2)  &  27  &  26  &  1262  &  4187  &  50  \\ \hline
(2)  &  29  &  28  &  1597  &  5251  &  162  \\ \hline
(2)  &  31  &  30  &  1986  &  6481  &  116  \\ \hline
(2)  &  37  &  36  &  3517  &  11287  &  210  \\ \hline
(2)  &  41  &  40  &  4881  &  15541  &  312  \\ \hline
(2)  &  43  &  42  &  5678  &  18019  &  246  \\ \hline
(2)  &  47  &  46  &  7522  &  23737  &  495  \\ \hline
(2)  &  49  &  48  &  8577  &  27001  &  188  \\ \hline
\end{array}$$

%%%%%%%%%%%%%%%%

$$\begin{array}{|c|c|c|c|c|c|}\hline
\mbox{Case}& q & n & \mbox{Min.~genus} & \mbox{Max.~genus} & \nu_i(q) \\ \hline
(3)  &  4  &  3  &  1  &  1  &  1  \\ \hline
(3)  &  8  &  7  &  8  &  64  &  6  \\ \hline
(3)  &  16  &  15  &  206  &  766  &  14  \\ \hline
(3)  &  32  &  31  &  2202  &  7162  &  90  \\ \hline
\end{array}$$

%%%%%%%%%%%%%%%%

$$\begin{array}{|c|c|c|c|c|c|}\hline
\mbox{Case}& q & n & \mbox{Min.~genus} & \mbox{Max.~genus} & \nu_i(q) \\ \hline
(4)  &  11  &  6  &  46  &  57  &  2  \\ \hline
(4)  &  13  &  7  &  15  &  155  &  9  \\ \hline
(4)  &  17  &  9  &  70  &  406  &  12  \\ \hline
(4)  &  19  &  10  &  116  &  515  &  8  \\ \hline
(4)  &  23  &  12  &  255  &  1014  &  10  \\ \hline
(4)  &  25  &  13  &  352  &  1352  &  18  \\ \hline
(4)  &  27  &  14  &  1406  &  1757  &  6  \\ \hline
(4)  &  29  &  15  &  611  &  2431  &  28  \\ \hline
(4)  &  31  &  16  &  777  &  2792  &  28  \\ \hline
(4)  &  37  &  19  &  1445  &  5321  &  81  \\ \hline
(4)  &  41  &  21  &  2052  &  7372  &  60  \\ \hline
(4)  &  43  &  22  &  2410  &  8129  &  50  \\ \hline
(4)  &  47  &  24  &  3245  &  10812  &  44  \\ \hline
(4)  &  49  &  25  &  3726  &  12350  &  60  \\ \hline
\end{array}$$

%%%%%%%%%%%%%%%%

$$\begin{array}{|c|c|c|c|c|c|}\hline
\mbox{Case}& q & n & \mbox{Min.~genus} & \mbox{Max.~genus} & \nu_i(q) \\ \hline
(5)  &  7  &  8  &  9  &  44  &  10  \\ \hline
(5)  &  9  &  10  &  26  &  110  &  7  \\ \hline
(5)  &  11  &  12  &  57  &  222  &  18  \\ \hline
(5)  &  13  &  14  &  106  &  392  &  33  \\ \hline
(5)  &  17  &  18  &  274  &  954  &  45  \\ \hline
(5)  &  19  &  20  &  401  &  1370  &  68  \\ \hline
(5)  &  23  &  24  &  761  &  2532  &  84  \\ \hline
(5)  &  25  &  26  &  1002  &  3302  &  69  \\ \hline
(5)  &  27  &  28  &  1289  &  4214  &  50  \\ \hline
(5)  &  29  &  30  &  1626  &  5280  &  108  \\ \hline
(5)  &  31  &  32  &  2017  &  6512  &  232  \\ \hline
(5)  &  37  &  38  &  3554  &  11324  &  315  \\ \hline
(5)  &  41  &  42  &  4922  &  15582  &  234  \\ \hline
(5)  &  43  &  44  &  5721  &  18062  &  410  \\ \hline
(5)  &  47  &  48  &  7569  &  23784  &  360  \\ \hline
(5)  &  49  &  50  &  8626  &  27050  &  235  \\ \hline
\end{array}$$

%%%%%%%%%%%%%%%%%%

$$\begin{array}{|c|c|c|c|c|c|}\hline
\mbox{Case}& q & n & \mbox{Min.~genus} & \mbox{Max.~genus} & \nu_i(q) \\ \hline
(6)  &  4  &  5  &  1  &  5  &  2  \\ \hline
(6)  &  8  &  9  &  16  &  72  &  6  \\ \hline
(6)  &  16  &  17  &  222  &  782  &  28  \\ \hline
(6)  &  32  &  33  &  2234  &  7194  &  60  \\ \hline 
\end{array}$$

$$\begin{array}{|c|c|c|c|c|c|}\hline
\mbox{Case}& q & n & \mbox{Min.~genus} & \mbox{Max.~genus} & \nu_i(q) \\ \hline
(7)  &  8  &  7  &  2290  &  4082  &  6  \\ \hline
(7)  &  32  &  31  &   4355842  &  7212802  &  90  \\ \hline
\end{array}$$

%%%%%%%%%%%%%%%%%%%%%

%\bigskip
\newpage

\section{Appendix B: data for maps with affine groups}

The first four tables illustrate the fact that, when the automorphism group $G$ is an affine group, the number $\phi(n)/d$ of Petrie dual pairs of maps ${\mathcal M}^{\pm}$ is unbounded for some infinite sequences of (odd) valencies $n$ (see the Example in Section~\ref{sec:affexs}).

\medskip

Mersenne numbers $n=2^e-1, e\ge 3$ (Case B).

$$\begin{array}{|r|r|r|r|}\hline
d & e & n & \phi(n)/d \\ \hline
           6&        3&        7&        1 \\ \hline
           8&        4&       15&        1 \\ \hline
          10&        5&       31&        3 \\ \hline
          12&        6&       63&        3 \\ \hline
          14&        7&      127&        9 \\ \hline
          16&        8&      255&        8 \\ \hline
          18&        9&      511&       24 \\ \hline
          20&       10&     1023&       30 \\ \hline
          22&       11&     2047&       88 \\ \hline
          24&       12&     4095&       72 \\ \hline
          26&       13&     8191&      315 \\ \hline
          28&       14&    16383&      378 \\ \hline
          30&       15&    32767&      900 \\ \hline
          32&       16&    65535&     1024 \\ \hline
          34&       17&   131071&     3855 \\ \hline
          36&       18&   262143&     3888 \\ \hline
          38&       19&   524287&    13797 \\ \hline
          40&       20&  1048575&    12000 \\ \hline
\end{array}$$

%%%%%%%%%%

\newpage

Mersenne numbers $n=2^e-1$, prime $e\ge 5$ (Case B).

$$\begin{array}{|r|r|r|r|}\hline
d & e & n & \phi(n)/d \\ \hline
                         10 &                         5 &                        31 &                         3 \\ \hline
                          14 &                         7 &                       127 &                         9 \\ \hline
                          22 &                        11 &                      2047 &                        88 \\ \hline
                          26 &                        13 &                      8191 &                       315 \\ \hline
                          34 &                        17 &                    131071 &                      3855 \\ \hline
                          38 &                        19 &                    524287 &                     13797 \\ \hline
                          46 &                        23 &                   8388607 &                    178480 \\ \hline
                          58 &                        29 &                 536870911 &                   9203904 \\ \hline
                          62 &                        31 &                2147483647 &                  34636833 \\ \hline
                          74 &                        37 &              137438953471 &                1848954528 \\ \hline
                          82 &                        41 &             2199023255551 &               26815350376 \\ \hline
                          86 &                        43 &             8796093022207 &              102032294580 \\ \hline
                          94 &                        47 &           140737488355327 &             1496238758400 \\ \hline
                         106 &                        53 &          9007199254740991 &            84958991520000 \\ \hline
                         118 &                        59 &        576460752303423487 &          4885233465012400 \\ \hline
                         122 &                        61 &       2305843009213693951 &         18900352534538475 \\ \hline
                         134 &                        67 &     147573952589676412927 &       1101298147967495880 \\ \hline
                         142 &                        71 &    2361183241434822606847 &      16627977798226714560 \\ \hline
                         146 &                        73 &    9444732965739290427391 &      64542566212975464960 \\ \hline
                         158 &                        79 &  604462909807314587353087 &    3824290813491610605444 \\ \hline
\end{array}$$

%%%%%%%%%%%%

\newpage

Valencies $n=2^e+1, e\ge 2$ (Case A).

$$\begin{array}{|r|r|r|r|}\hline
d & e & n & \phi(n)/d \\ \hline
         4 &        2 &        5 &        1 \\ \hline
          6 &        3 &        9 &        1 \\ \hline
          8 &        4 &       17 &        2 \\ \hline
         10 &        5 &       33 &        2 \\ \hline
         12 &        6 &       65 &        4 \\ \hline
         14 &        7 &      129 &        6 \\ \hline
         16 &        8 &      257 &       16 \\ \hline
         18 &        9 &      513 &       18 \\ \hline
         20 &       10 &     1025 &       40 \\ \hline
         22 &       11 &     2049 &       62 \\ \hline
         24 &       12 &     4097 &      160 \\ \hline
         26 &       13 &     8193 &      210 \\ \hline
         28 &       14 &    16385 &      448 \\ \hline
         30 &       15 &    32769 &      660 \\ \hline
         32 &       16 &    65537 &     2048 \\ \hline
         34 &       17 &   131073 &     2570 \\ \hline
         36 &       18 &   262145 &     5184 \\ \hline
         38 &       19 &   524289 &     9198 \\ \hline
         40 &       20 &  1048577 &    24672 \\ \hline
\end{array}$$

%%%%%%%%%%%%

\newpage

Valencies $n=2^e+1$, prime $e\ge 3$ (Case A).

$$\begin{array}{|r|r|r|r|}\hline
d & e & n & \phi(n)/d \\ \hline
                          6&                       3&                       9&                       1 \\ \hline
                         10&                       5&                      33&                       2 \\ \hline
                         14&                       7&                     129&                       6 \\ \hline
                         22&                      11&                    2049&                      62 \\ \hline
                         26&                      13&                    8193&                     210 \\ \hline
                         34&                      17&                  131073&                    2570 \\ \hline
                         38&                      19&                  524289&                    9198 \\ \hline
                         46&                      23&                 8388609&                  121574 \\ \hline
                         58&                      29&               536870913&                 6066336 \\ \hline
                         62&                      31&              2147483649&                23091222 \\ \hline
                         74&                      37&            137438953473&              1237491936 \\ \hline
                         82&                      41&           2199023255553&             17662837392 \\ \hline
                         86&                      43&           8796093022209&             68186767614 \\ \hline
                         94&                      47&         140737488355329&            994611225120 \\ \hline
                        106&                      53&        9007199254740993&          56119621524864 \\ \hline
                        118&                      59&      576460752303423489&        3255603211526400 \\ \hline
                        122&                      61&     2305843009213693953&       12600235023025650 \\ \hline
                        134&                      67&   147573952589676412929&      734198668907188464 \\ \hline
                        142&                      71&  2361183241434822606849&    11085367134163656960 \\ \hline
                        146&                      73&  9444732965739290427393&    43102032929785695024 \\ \hline
\end{array}$$

%%%%%%%%%%%%

\newpage

The next two tables (the first in two parts, on account of its length) show data concerning the maps ${\mathcal M}^{\pm}$ for admissible valencies $n \leq 50$ with $p=2$ and $p=3$. In the final column, $\phi(n)/d$ is the number of Petrie dual pairs of maps.

$$\begin{array}{|r|r|r|r|r|r|r|r|r|r|}\hline
p & d & e & n & \mbox{case} & \mbox{map} & \mbox{type} & \mbox{genus}  & \mbox{self-dual} & \phi(n)/d \\ \hline
2 & 2 & 1 & 3 & A & {\mathcal M}^+ & \{ 3 , 3 \}_{4} & 0 & \mbox{yes} & 1\\ \hline
 &  &  &  &  & {\mathcal M}^- &\{ 4 , 3 \}_{3} & 1 & & \\ \hline
 2 & 4 & 2 & 5 & A & {\mathcal M}^+ & \{ 5 , 5 \}_{4} & 5 & \mbox{yes} & 1\\ \hline
 &  &  &  &  & {\mathcal M}^- &\{ 4 , 5 \}_{5} & 6 & & \\ \hline
2 & 6 & 3 & 7 & B & {\mathcal M}^+ & \{ 7 , 7 \}_{4} & 49 & \mbox{yes} & 1\\ \hline
 &  &  &  &  & {\mathcal M}^- &\{ 4 , 7 \}_{7} & 50 & & \\ \hline
2 & 6 & 3 & 9 & A & {\mathcal M}^+ & \{ 9 , 9 \}_{4} & 81 & \mbox{yes} & 1\\ \hline
 &  &  &  &  & {\mathcal M}^- &\{ 4 , 9 \}_{9} & 82 & & \\ \hline
2 & 10 & 5 & 11 & A & {\mathcal M}^+ & \{ 11 , 11 \}_{4} & 1793 & \mbox{yes} & 1\\ \hline
 &  &  &  &  & {\mathcal M}^- &\{ 4 , 11 \}_{11} & 1794 & & \\ \hline
2 & 12 & 6 & 13 & A & {\mathcal M}^+ & \{ 13 , 13 \}_{4} & 9217 & \mbox{yes} & 1\\ \hline
 &  &  &  &  & {\mathcal M}^- &\{ 4 , 13 \}_{13} & 9218 & & \\ \hline
2 & 8 & 4 & 15 & B & {\mathcal M}^+ & \{ 15 , 15 \}_{4} & 705 & \mbox{yes} & 1\\ \hline
 &  &  &  &  & {\mathcal M}^- &\{ 4 , 15 \}_{15} & 706 & & \\ \hline
2 & 8 & 4 & 17 & A & {\mathcal M}^+ & \{ 17 , 17 \}_{4} & 833 & \mbox{yes} & 2\\ \hline
 &  &  &  &  & {\mathcal M}^- &\{ 4 , 17 \}_{17} & 834 & & \\ \hline
2 & 18 & 9 & 19 & A & {\mathcal M}^+ & \{ 19 , 19 \}_{4} & 983041 & \mbox{yes} & 1\\ \hline
 &  &  &  &  & {\mathcal M}^- &\{ 4 , 19 \}_{19} & 983042 & & \\ \hline
2 & 12 & 6 & 21 & B & {\mathcal M}^+ & \{ 21 , 21 \}_{4} & 17409 & \mbox{yes} & 1\\ \hline
 &  &  &  &  & {\mathcal M}^- &\{ 4 , 21 \}_{21} & 17410 & & \\ \hline
2 & 22 & 11 & 23 & B & {\mathcal M}^+ & \{ 23 , 23 \}_{4} & 19922945 & \mbox{yes} & 1\\ \hline
 &  &  &  &  & {\mathcal M}^- &\{ 4 , 23 \}_{23} & 19922946 & & \\ \hline
2 & 20 & 10 & 25 & A & {\mathcal M}^+ & \{ 25 , 25 \}_{4} & 5505025 & \mbox{yes} & 1\\ \hline
 &  &  &  &  & {\mathcal M}^- &\{ 4 , 25 \}_{25} & 5505026 & & \\ \hline
2 & 18 & 9 & 27 & A & {\mathcal M}^+ & \{ 27 , 27 \}_{4} & 1507329 & \mbox{yes} & 1\\ \hline
 &  &  &  &  & {\mathcal M}^- &\{ 4 , 27 \}_{27} & 1507330 & & \\ \hline
2 & 28 & 14 & 29 & A & {\mathcal M}^+ & \{ 29 , 29 \}_{4} & 1677721601 & \mbox{yes} & 1\\ \hline
 &  &  &  &  & {\mathcal M}^- &\{ 4 , 29 \}_{29} & 1677721602 & & \\ \hline
2 & 10 & 5 & 31 & B & {\mathcal M}^+ & \{ 31 , 31 \}_{4} & 6913 & \mbox{yes} & 3\\ \hline
 &  &  &  &  & {\mathcal M}^- &\{ 4 , 31 \}_{31} & 6914 & & \\ \hline
2 & 10 & 5 & 33 & A & {\mathcal M}^+ & \{ 33 , 33 \}_{4} & 7425 & \mbox{yes} & 2\\ \hline
 &  &  &  &  & {\mathcal M}^- &\{ 4 , 33 \}_{33} & 7426 & & \\ \hline
2 & 24 & 12 & 35 & B & {\mathcal M}^+ & \{ 35 , 35 \}_{4} & 130023425 & \mbox{yes} & 1\\ \hline
 &  &  &  &  & {\mathcal M}^- &\{ 4 , 35 \}_{35} & 130023426 & & \\ \hline
2 & 36 & 18 & 37 & A & {\mathcal M}^+ & \{ 37 , 37 \}_{4} & 566935683073 & \mbox{yes} & 1\\ \hline

\end{array}$$

$$\begin{array}{|r|r|r|r|r|r|r|r|r|r|}\hline
p & d & e & n & \mbox{case} & \mbox{map} & \mbox{type} & \mbox{genus}  & \mbox{self-dual} & \phi(n)/d \\ \hline
 &  &  &  &  & {\mathcal M}^- &\{ 4 , 37 \}_{37} & 566935683074 & & \\ \hline
2 & 24 & 12 & 39 & B & {\mathcal M}^+ & \{ 39 , 39 \}_{4} & 146800641 & \mbox{yes} & 1\\ \hline
 &  &  &  &  & {\mathcal M}^- &\{ 4 , 39 \}_{39} & 146800642 & & \\ \hline
2 & 20 & 10 & 41 & A & {\mathcal M}^+ & \{ 41 , 41 \}_{4} & 9699329 & \mbox{yes} & 2\\ \hline
 &  &  &  &  & {\mathcal M}^- &\{ 4 , 41 \}_{41} & 9699330 & & \\ \hline
2 & 14 & 7 & 43 & A & {\mathcal M}^+ & \{ 43 , 43 \}_{4} & 159745 & \mbox{yes} & 3\\ \hline
 &  &  &  &  & {\mathcal M}^- &\{ 4 , 43 \}_{43} & 159746 & & \\ \hline
2 & 24 & 12 & 45 & B & {\mathcal M}^+ & \{ 45 , 45 \}_{4} & 171966465 & \mbox{yes} & 1\\ \hline
 &  &  &  &  & {\mathcal M}^- &\{ 4 , 45 \}_{45} & 171966466 & & \\ \hline
2 & 46 & 23 & 47 & B & {\mathcal M}^+ & \{ 47 , 47 \}_{4} & 756463999909889 & \mbox{yes} & 1\\ \hline
 &  &  &  &  & {\mathcal M}^- &\{ 4 , 47 \}_{47} & 756463999909890 & & \\ \hline
2 & 42 & 21 & 49 & B & {\mathcal M}^+ & \{ 49 , 49 \}_{4} & 49478023249921 & \mbox{yes} & 1\\ \hline
 &  &  &  &  & {\mathcal M}^- &\{ 4 , 49 \}_{49} & 49478023249922 & & \\ \hline
\end{array}$$

\newpage
$$\begin{array}{|r|r|r|r|r|r|r|r|r|r|}\hline
p & d & e & n & \mbox{case} & \mbox{map} & \mbox{type} & \mbox{genus}  & \mbox{self-dual} & \phi(n)/d \\ \hline
3 & 2 & 1 & 4 & A & {\mathcal M}^+ & \{ 4 , 4 \}_{6} & 1 & \mbox{yes} & 1\\ \hline
 &  &  &  &  & {\mathcal M}^- &\{ 6 , 4 \}_{4} & 5 & & \\ \hline
3 & 4 & 2 & 8 & B & {\mathcal M}^+ & \{ 8 , 8 \}_{6} & 82 & \mbox{yes} & 1\\ \hline
 &  &  &  &  & {\mathcal M}^- &\{ 6 , 8 \}_{8} & 137 & & \\ \hline
3 & 4 & 2 & 10 & A & {\mathcal M}^+ & \{ 5 , 10 \}_{6} & 82 & \mbox{no} & 1\\ \hline
 &  &  &  &  & {\mathcal M}^- &\{ 6 , 10 \}_{5} & 191 & & \\ \hline
3 & 6 & 3 & 14 & A & {\mathcal M}^+ & \{ 7 , 14 \}_{6} & 1459 & \mbox{no} & 1\\ \hline
 &  &  &  &  & {\mathcal M}^- &\{ 6 , 14 \}_{7} & 2675 & & \\ \hline
3 & 8 & 4 & 16 & B & {\mathcal M}^+ & \{ 16 , 16 \}_{6} & 19684 & \mbox{yes} & 1\\ \hline
 &  &  &  &  & {\mathcal M}^- &\{ 6 , 16 \}_{16} & 28433 & & \\ \hline
3 & 8 & 4 & 20 & B & {\mathcal M}^+ & \{ 20 , 20 \}_{6} & 26245 & \mbox{yes} & 1\\ \hline
 &  &  &  &  & {\mathcal M}^- &\{ 6 , 20 \}_{20} & 37181 & & \\ \hline
3 & 10 & 5 & 22 & B & {\mathcal M}^+ & \{ 11 , 22 \}_{6} & 236197 & \mbox{no} & 1\\ \hline
 &  &  &  &  & {\mathcal M}^- &\{ 6 , 22 \}_{11} & 373979 & & \\ \hline
3 & 6 & 3 & 26 & B & {\mathcal M}^+ & \{ 13 , 26 \}_{6} & 3646 & \mbox{no} & 2\\ \hline
 &  &  &  &  & {\mathcal M}^- &\{ 6 , 26 \}_{13} & 5591 & & \\ \hline
3 & 6 & 3 & 28 & A & {\mathcal M}^+ & \{ 28 , 28 \}_{6} & 4375 & \mbox{no} & 2\\ \hline
 &  &  &  &  & {\mathcal M}^- &\{ 6 , 28 \}_{28} & 6077 & & \\ \hline
3 & 16 & 8 & 32 & B & {\mathcal M}^+ & \{ 32 , 32 \}_{6} & 301327048 & \mbox{yes} & 1\\ \hline
 &  &  &  &  & {\mathcal M}^- &\{ 6 , 32 \}_{32} & 416118305 & & \\ \hline
3 & 16 & 8 & 34 & A & {\mathcal M}^+ & \{ 17 , 34 \}_{6} & 301327048 & \mbox{no} & 1\\ \hline
 &  &  &  &  & {\mathcal M}^- &\{ 6 , 34 \}_{17} & 444816119 & & \\ \hline
3 & 18 & 9 & 38 & A & {\mathcal M}^+ & \{ 19 , 38 \}_{6} & 3099363913 & \mbox{no} & 1\\ \hline
 &  &  &  &  & {\mathcal M}^- &\{ 6 , 38 \}_{19} & 4519905707 & & \\ \hline
3 & 8 & 4 & 40 & B & {\mathcal M}^+ & \{ 40 , 40 \}_{6} & 59050 & \mbox{no} & 2\\ \hline
 &  &  &  &  & {\mathcal M}^- &\{ 6 , 40 \}_{40} & 80921 & & \\ \hline
3 & 20 & 10 & 44 & B & {\mathcal M}^+ & \{ 44 , 44 \}_{6} & 34867844011 & \mbox{yes} & 1\\ \hline
 &  &  &  &  & {\mathcal M}^- &\{ 6 , 44 \}_{44} & 47652720149 & & \\ \hline
3 & 22 & 11 & 46 & B & {\mathcal M}^+ & \{ 23 , 46 \}_{6} & 313810596091 & \mbox{no} & 1\\ \hline
 &  &  &  &  & {\mathcal M}^- &\{ 6 , 46 \}_{23} & 449795187731 & & \\ \hline
3 & 20 & 10 & 50 & A & {\mathcal M}^+ & \{ 25 , 50 \}_{6} & 38354628412 & \mbox{no} & 1\\ \hline
 &  &  &  &  & {\mathcal M}^- &\{ 6 , 50 \}_{25} & 54626288951 & & \\ \hline
\end{array}$$

%%%%%%%%%%%%

\newpage

The following tables show properties of maps of a given even valency $n=4,\ldots, 30$. The expressions for the type and genus, as functions of the prime $p$, depend only on $n$ and $p$ mod~$n$, so by Dirichlet's Theorem each row represents an infinite sequence of maps.

$$\begin{array}{|r|r|r|r|r|r|r|r|r|r|}\hline
p\; \mbox{mod}\; n & d & e & n & \mbox{case} & \mbox{map} & \mbox{type} & \mbox{genus}  & \mbox{self-dual} & \phi(n)/d \\ \hline
1 & 2 & 1 & 4 & B & {\mathcal M}^+ & \{ 4 , 4 \}_{2p} & 1 & \mbox{yes} & 1\\ \hline
 &  &  &  &  & {\mathcal M}^- &\{ 2p , 4 \}_{4} & 2+ p^2( 2p- 4)/2p & & \\ \hline
3 & 2 & 1 & 4 & A & {\mathcal M}^+ & \{ 4 , 4 \}_{2p} & 1 & \mbox{yes} & 1\\ \hline
 &  &  &  &  & {\mathcal M}^- &\{ 2p , 4 \}_{4} & 2+ p^2( 2p- 4)/2p & & \\ \hline
\end{array}$$

$$\begin{array}{|r|r|r|r|r|r|r|r|r|r|}\hline
p\; \mbox{mod}\; n & d & e & n & \mbox{case} & \mbox{map} & \mbox{type} & \mbox{genus}  & \mbox{self-dual} & \phi(n)/d \\ \hline
1 & 2 & 1 & 6 & B & {\mathcal M}^+ & \{ 3 , 6 \}_{2p} & 1 & \mbox{no} & 1\\ \hline
 &  &  &  &  & {\mathcal M}^- &\{ 2p , 6 \}_{3} & 2+ p^{2}( 4p- 6)/2p & & \\ \hline
5 & 2 & 1 & 6 & A & {\mathcal M}^+ & \{ 3 , 6 \}_{2p} & 1 & \mbox{no} & 1\\ \hline
 &  &  &  &  & {\mathcal M}^- &\{ 2p , 6 \}_{3} & 2+ p^{2}( 4p- 6)/2p & & \\ \hline
\end{array}$$

$$\begin{array}{|r|r|r|r|r|r|r|r|r|r|}\hline
p\; \mbox{mod}\; n & d & e & n & \mbox{case} & \mbox{map} & \mbox{type} & \mbox{genus}  & \mbox{self-dual} & \phi(n)/d \\ \hline
1 & 2 & 1 & 8 & B & {\mathcal M}^+ & \{ 8 , 8 \}_{2p} & 1+p^{2} & \mbox{no} & 2\\ \hline
 &  &  &  &  & {\mathcal M}^- &\{ 2p , 8 \}_{8} & 2+ p^{2}( 6p- 8)/2p & & \\ \hline
3 & 4 & 2 & 8 & B & {\mathcal M}^+ & \{ 8 , 8 \}_{2p} & 1+p^{4} & \mbox{yes} & 1\\ \hline
 &  &  &  &  & {\mathcal M}^- &\{ 2p , 8 \}_{8} & 2+ p^{4}( 6p- 8)/2p & & \\ \hline
5 & 4 & 2 & 8 & B & {\mathcal M}^+ & \{ 8 , 8 \}_{2p} & 1+p^{4} & \mbox{yes} & 1\\ \hline
 &  &  &  &  & {\mathcal M}^- &\{ 2p , 8 \}_{8} & 2+ p^{4}( 6p- 8)/2p & & \\ \hline
7 & 2 & 1 & 8 & A & {\mathcal M}^+ & \{ 8 , 8 \}_{2p} & 1+p^{2} & \mbox{no} & 2\\ \hline
 &  &  &  &  & {\mathcal M}^- &\{ 2p , 8 \}_{8} & 2+ p^{2}( 6p- 8)/2p & & \\ \hline
\end{array}$$

$$\begin{array}{|r|r|r|r|r|r|r|r|r|r|}\hline
p\; \mbox{mod}\; n & d & e & n & \mbox{case} & \mbox{map} & \mbox{type} & \mbox{genus}  & \mbox{self-dual} & \phi(n)/d \\ \hline
1 & 2 & 1 & 10 & B & {\mathcal M}^+ & \{ 5 , 10 \}_{2p} & 1+p^{2} & \mbox{no} & 2\\ \hline
 &  &  &  &  & {\mathcal M}^- &\{ 2p , 10 \}_{5} & 2+ p^{2}( 8p- 10)/2p & & \\ \hline
3 & 4 & 2 & 10 & A & {\mathcal M}^+ & \{ 5 , 10 \}_{2p} & 1+p^{4} & \mbox{no} & 1\\ \hline
 &  &  &  &  & {\mathcal M}^- &\{ 2p , 10 \}_{5} & 2+ p^{4}( 8p- 10)/2p & & \\ \hline
7 & 4 & 2 & 10 & A & {\mathcal M}^+ & \{ 5 , 10 \}_{2p} & 1+p^{4} & \mbox{no} & 1\\ \hline
 &  &  &  &  & {\mathcal M}^- &\{ 2p , 10 \}_{5} & 2+ p^{4}( 8p- 10)/2p & & \\ \hline
9 & 2 & 1 & 10 & A & {\mathcal M}^+ & \{ 5 , 10 \}_{2p} & 1+p^{2} & \mbox{no} & 2\\ \hline
 &  &  &  &  & {\mathcal M}^- &\{ 2p , 10 \}_{5} & 2+ p^{2}( 8p- 10)/2p & & \\ \hline
\end{array}$$

$$\begin{array}{|r|r|r|r|r|r|r|r|r|r|}\hline
p\; \mbox{mod}\; n & d & e & n & \mbox{case} & \mbox{map} & \mbox{type} & \mbox{genus}  & \mbox{self-dual} & \phi(n)/d \\ \hline
1 & 2 & 1 & 12 & B & {\mathcal M}^+ & \{ 12 , 12 \}_{2p} & 1+2p^{2} & \mbox{no} & 2\\ \hline
 &  &  &  &  & {\mathcal M}^- &\{ 2p , 12 \}_{12} & 2+ p^{2}( 10p- 12)/2p & & \\ \hline
5 & 4 & 2 & 12 & B & {\mathcal M}^+ & \{ 12 , 12 \}_{2p} & 1+2p^{4} & \mbox{yes} & 1\\ \hline
 &  &  &  &  & {\mathcal M}^- &\{ 2p , 12 \}_{12} & 2+ p^{4}( 10p- 12)/2p & & \\ \hline
7 & 4 & 2 & 12 & B & {\mathcal M}^+ & \{ 12 , 12 \}_{2p} & 1+2p^{4} & \mbox{yes} & 1\\ \hline
 &  &  &  &  & {\mathcal M}^- &\{ 2p , 12 \}_{12} & 2+ p^{4}( 10p- 12)/2p & & \\ \hline
11 & 2 & 1 & 12 & A & {\mathcal M}^+ & \{ 12 , 12 \}_{2p} & 1+2p^{2} & \mbox{no} & 2\\ \hline
 &  &  &  &  & {\mathcal M}^- &\{ 2p , 12 \}_{12} & 2+ p^{2}( 10p- 12)/2p & & \\ \hline
\end{array}$$

$$\begin{array}{|r|r|r|r|r|r|r|r|r|r|}\hline
p\; \mbox{mod}\; n & d & e & n & \mbox{case} & \mbox{map} & \mbox{type} & \mbox{genus}  & \mbox{self-dual} & \phi(n)/d \\ \hline
1 & 2 & 1 & 14 & B & {\mathcal M}^+ & \{ 7 , 14 \}_{2p} & 1+2p^{2} & \mbox{no} & 3\\ \hline
 &  &  &  &  & {\mathcal M}^- &\{ 2p , 14 \}_{7} & 2+ p^{2}( 12p- 14)/2p & & \\ \hline
3 & 6 & 3 & 14 & A & {\mathcal M}^+ & \{ 7 , 14 \}_{2p} & 1+2p^{6} & \mbox{no} & 1\\ \hline
 &  &  &  &  & {\mathcal M}^- &\{ 2p , 14 \}_{7} & 2+ p^{6}( 12p- 14)/2p & & \\ \hline
5 & 6 & 3 & 14 & A & {\mathcal M}^+ & \{ 7 , 14 \}_{2p} & 1+2p^{6} & \mbox{no} & 1\\ \hline
 &  &  &  &  & {\mathcal M}^- &\{ 2p , 14 \}_{7} & 2+ p^{6}( 12p- 14)/2p & & \\ \hline
9 & 6 & 3 & 14 & B & {\mathcal M}^+ & \{ 7 , 14 \}_{2p} & 1+2p^{6} & \mbox{no} & 1\\ \hline
 &  &  &  &  & {\mathcal M}^- &\{ 2p , 14 \}_{7} & 2+ p^{6}( 12p- 14)/2p & & \\ \hline
11 & 6 & 3 & 14 & B & {\mathcal M}^+ & \{ 7 , 14 \}_{2p} & 1+2p^{6} & \mbox{no} & 1\\ \hline
 &  &  &  &  & {\mathcal M}^- &\{ 2p , 14 \}_{7} & 2+ p^{6}( 12p- 14)/2p & & \\ \hline
13 & 2 & 1 & 14 & A & {\mathcal M}^+ & \{ 7 , 14 \}_{2p} & 1+2p^{2} & \mbox{no} & 3\\ \hline
 &  &  &  &  & {\mathcal M}^- &\{ 2p , 14 \}_{7} & 2+ p^{2}( 12p- 14)/2p & & \\ \hline
\end{array}$$

$$\begin{array}{|r|r|r|r|r|r|r|r|r|r|}\hline
p\; \mbox{mod}\; n & d & e & n & \mbox{case} & \mbox{map} & \mbox{type} & \mbox{genus}  & \mbox{self-dual} & \phi(n)/d \\ \hline
1 & 2 & 1 & 16 & B & {\mathcal M}^+ & \{ 16 , 16 \}_{2p} & 1+3p^{2} & \mbox{no} & 4\\ \hline
 &  &  &  &  & {\mathcal M}^- &\{ 2p , 16 \}_{16} & 2+ p^{2}( 14p- 16)/2p & & \\ \hline
3 & 8 & 4 & 16 & B & {\mathcal M}^+ & \{ 16 , 16 \}_{2p} & 1+3p^{8} & \mbox{yes} & 1\\ \hline
 &  &  &  &  & {\mathcal M}^- &\{ 2p , 16 \}_{16} & 2+ p^{8}( 14p- 16)/2p & & \\ \hline
5 & 8 & 4 & 16 & B & {\mathcal M}^+ & \{ 16 , 16 \}_{2p} & 1+3p^{8} & \mbox{yes} & 1\\ \hline
 &  &  &  &  & {\mathcal M}^- &\{ 2p , 16 \}_{16} & 2+ p^{8}( 14p- 16)/2p & & \\ \hline
7 & 4 & 2 & 16 & B & {\mathcal M}^+ & \{ 16 , 16 \}_{2p} & 1+3p^{4} & \mbox{yes} & 2\\ \hline
 &  &  &  &  & {\mathcal M}^- &\{ 2p , 16 \}_{16} & 2+ p^{4}( 14p- 16)/2p & & \\ \hline
9 & 4 & 2 & 16 & B & {\mathcal M}^+ & \{ 16 , 16 \}_{2p} & 1+3p^{4} & \mbox{yes} & 2\\ \hline
 &  &  &  &  & {\mathcal M}^- &\{ 2p , 16 \}_{16} & 2+ p^{4}( 14p- 16)/2p & & \\ \hline
11 & 8 & 4 & 16 & B & {\mathcal M}^+ & \{ 16 , 16 \}_{2p} & 1+3p^{8} & \mbox{yes} & 1\\ \hline
 &  &  &  &  & {\mathcal M}^- &\{ 2p , 16 \}_{16} & 2+ p^{8}( 14p- 16)/2p & & \\ \hline
13 & 8 & 4 & 16 & B & {\mathcal M}^+ & \{ 16 , 16 \}_{2p} & 1+3p^{8} & \mbox{yes} & 1\\ \hline
 &  &  &  &  & {\mathcal M}^- &\{ 2p , 16 \}_{16} & 2+ p^{8}( 14p- 16)/2p & & \\ \hline
15 & 2 & 1 & 16 & A & {\mathcal M}^+ & \{ 16 , 16 \}_{2p} & 1+3p^{2} & \mbox{no} & 4\\ \hline
 &  &  &  &  & {\mathcal M}^- &\{ 2p , 16 \}_{16} & 2+ p^{2}( 14p- 16)/2p & & \\ \hline
\end{array}$$

$$\begin{array}{|r|r|r|r|r|r|r|r|r|r|}\hline
p\; \mbox{mod}\; n & d & e & n & \mbox{case} & \mbox{map} & \mbox{type} & \mbox{genus}  & \mbox{self-dual} & \phi(n)/d \\ \hline
1 & 2 & 1 & 18 & B & {\mathcal M}^+ & \{ 9 , 18 \}_{2p} & 1+3p^{2} & \mbox{no} & 3\\ \hline
 &  &  &  &  & {\mathcal M}^- &\{ 2p , 18 \}_{9} & 2+ p^{2}( 16p- 18)/2p & & \\ \hline
5 & 6 & 3 & 18 & A & {\mathcal M}^+ & \{ 9 , 18 \}_{2p} & 1+3p^{6} & \mbox{no} & 1\\ \hline
 &  &  &  &  & {\mathcal M}^- &\{ 2p , 18 \}_{9} & 2+ p^{6}( 16p- 18)/2p & & \\ \hline
7 & 6 & 3 & 18 & B & {\mathcal M}^+ & \{ 9 , 18 \}_{2p} & 1+3p^{6} & \mbox{no} & 1\\ \hline
 &  &  &  &  & {\mathcal M}^- &\{ 2p , 18 \}_{9} & 2+ p^{6}( 16p- 18)/2p & & \\ \hline
11 & 6 & 3 & 18 & A & {\mathcal M}^+ & \{ 9 , 18 \}_{2p} & 1+3p^{6} & \mbox{no} & 1\\ \hline
 &  &  &  &  & {\mathcal M}^- &\{ 2p , 18 \}_{9} & 2+ p^{6}( 16p- 18)/2p & & \\ \hline
13 & 6 & 3 & 18 & B & {\mathcal M}^+ & \{ 9 , 18 \}_{2p} & 1+3p^{6} & \mbox{no} & 1\\ \hline
 &  &  &  &  & {\mathcal M}^- &\{ 2p , 18 \}_{9} & 2+ p^{6}( 16p- 18)/2p & & \\ \hline
17 & 2 & 1 & 18 & A & {\mathcal M}^+ & \{ 9 , 18 \}_{2p} & 1+3p^{2} & \mbox{no} & 3\\ \hline
 &  &  &  &  & {\mathcal M}^- &\{ 2p , 18 \}_{9} & 2+ p^{2}( 16p- 18)/2p & & \\ \hline
\end{array}$$

$$\begin{array}{|r|r|r|r|r|r|r|r|r|r|}\hline
p\; \mbox{mod}\; n & d & e & n & \mbox{case} & \mbox{map} & \mbox{type} & \mbox{genus}  & \mbox{self-dual} & \phi(n)/d \\ \hline
1 & 2 & 1 & 20 & B & {\mathcal M}^+ & \{ 20 , 20 \}_{2p} & 1+4p^{2} & \mbox{no} & 4\\ \hline
 &  &  &  &  & {\mathcal M}^- &\{ 2p , 20 \}_{20} & 2+ p^{2}( 18p- 20)/2p & & \\ \hline
3 & 8 & 4 & 20 & B & {\mathcal M}^+ & \{ 20 , 20 \}_{2p} & 1+4p^{8} & \mbox{yes} & 1\\ \hline
 &  &  &  &  & {\mathcal M}^- &\{ 2p , 20 \}_{20} & 2+ p^{8}( 18p- 20)/2p & & \\ \hline
7 & 8 & 4 & 20 & B & {\mathcal M}^+ & \{ 20 , 20 \}_{2p} & 1+4p^{8} & \mbox{yes} & 1\\ \hline
 &  &  &  &  & {\mathcal M}^- &\{ 2p , 20 \}_{20} & 2+ p^{8}( 18p- 20)/2p & & \\ \hline
9 & 4 & 2 & 20 & B & {\mathcal M}^+ & \{ 20 , 20 \}_{2p} & 1+4p^{4} & \mbox{yes} & 2\\ \hline
 &  &  &  &  & {\mathcal M}^- &\{ 2p , 20 \}_{20} & 2+ p^{4}( 18p- 20)/2p & & \\ \hline
11 & 4 & 2 & 20 & B & {\mathcal M}^+ & \{ 20 , 20 \}_{2p} & 1+4p^{4} & \mbox{yes} & 2\\ \hline
 &  &  &  &  & {\mathcal M}^- &\{ 2p , 20 \}_{20} & 2+ p^{4}( 18p- 20)/2p & & \\ \hline
13 & 8 & 4 & 20 & B & {\mathcal M}^+ & \{ 20 , 20 \}_{2p} & 1+4p^{8} & \mbox{yes} & 1\\ \hline
 &  &  &  &  & {\mathcal M}^- &\{ 2p , 20 \}_{20} & 2+ p^{8}( 18p- 20)/2p & & \\ \hline
17 & 8 & 4 & 20 & B & {\mathcal M}^+ & \{ 20 , 20 \}_{2p} & 1+4p^{8} & \mbox{yes} & 1\\ \hline
 &  &  &  &  & {\mathcal M}^- &\{ 2p , 20 \}_{20} & 2+ p^{8}( 18p- 20)/2p & & \\ \hline
19 & 2 & 1 & 20 & A & {\mathcal M}^+ & \{ 20 , 20 \}_{2p} & 1+4p^{2} & \mbox{no} & 4\\ \hline
 &  &  &  &  & {\mathcal M}^- &\{ 2p , 20 \}_{20} & 2+ p^{2}( 18p- 20)/2p & & \\ \hline
\end{array}$$

$$\begin{array}{|r|r|r|r|r|r|r|r|r|r|}\hline
p\; \mbox{mod}\; n & d & e & n & \mbox{case} & \mbox{map} & \mbox{type} & \mbox{genus}  & \mbox{self-dual} & \phi(n)/d \\ \hline
1 & 2 & 1 & 22 & B & {\mathcal M}^+ & \{ 11 , 22 \}_{2p} & 1+4p^{2} & \mbox{no} & 5\\ \hline
 &  &  &  &  & {\mathcal M}^- &\{ 2p , 22 \}_{11} & 2+ p^{2}( 20p- 22)/2p & & \\ \hline
3 & 10 & 5 & 22 & B & {\mathcal M}^+ & \{ 11 , 22 \}_{2p} & 1+4p^{10} & \mbox{no} & 1\\ \hline
 &  &  &  &  & {\mathcal M}^- &\{ 2p , 22 \}_{11} & 2+ p^{10}( 20p- 22)/2p & & \\ \hline
5 & 10 & 5 & 22 & B & {\mathcal M}^+ & \{ 11 , 22 \}_{2p} & 1+4p^{10} & \mbox{no} & 1\\ \hline
 &  &  &  &  & {\mathcal M}^- &\{ 2p , 22 \}_{11} & 2+ p^{10}( 20p- 22)/2p & & \\ \hline
7 & 10 & 5 & 22 & A & {\mathcal M}^+ & \{ 11 , 22 \}_{2p} & 1+4p^{10} & \mbox{no} & 1\\ \hline
 &  &  &  &  & {\mathcal M}^- &\{ 2p , 22 \}_{11} & 2+ p^{10}( 20p- 22)/2p & & \\ \hline
9 & 10 & 5 & 22 & B & {\mathcal M}^+ & \{ 11 , 22 \}_{2p} & 1+4p^{10} & \mbox{no} & 1\\ \hline
 &  &  &  &  & {\mathcal M}^- &\{ 2p , 22 \}_{11} & 2+ p^{10}( 20p- 22)/2p & & \\ \hline
13 & 10 & 5 & 22 & A & {\mathcal M}^+ & \{ 11 , 22 \}_{2p} & 1+4p^{10} & \mbox{no} & 1\\ \hline
 &  &  &  &  & {\mathcal M}^- &\{ 2p , 22 \}_{11} & 2+ p^{10}( 20p- 22)/2p & & \\ \hline
15 & 10 & 5 & 22 & B & {\mathcal M}^+ & \{ 11 , 22 \}_{2p} & 1+4p^{10} & \mbox{no} & 1\\ \hline
 &  &  &  &  & {\mathcal M}^- &\{ 2p , 22 \}_{11} & 2+ p^{10}( 20p- 22)/2p & & \\ \hline
17 & 10 & 5 & 22 & A & {\mathcal M}^+ & \{ 11 , 22 \}_{2p} & 1+4p^{10} & \mbox{no} & 1\\ \hline
 &  &  &  &  & {\mathcal M}^- &\{ 2p , 22 \}_{11} & 2+ p^{10}( 20p- 22)/2p & & \\ \hline
19 & 10 & 5 & 22 & A & {\mathcal M}^+ & \{ 11 , 22 \}_{2p} & 1+4p^{10} & \mbox{no} & 1\\ \hline
 &  &  &  &  & {\mathcal M}^- &\{ 2p , 22 \}_{11} & 2+ p^{10}( 20p- 22)/2p & & \\ \hline
21 & 2 & 1 & 22 & A & {\mathcal M}^+ & \{ 11 , 22 \}_{2p} & 1+4p^{2} & \mbox{no} & 5\\ \hline
 &  &  &  &  & {\mathcal M}^- &\{ 2p , 22 \}_{11} & 2+ p^{2}( 20p- 22)/2p & & \\ \hline
\end{array}$$

$$\begin{array}{|r|r|r|r|r|r|r|r|r|r|}\hline
p\; \mbox{mod}\; n & d & e & n & \mbox{case} & \mbox{map} & \mbox{type} & \mbox{genus}  & \mbox{self-dual} & \phi(n)/d \\ \hline
1 & 2 & 1 & 24 & B & {\mathcal M}^+ & \{ 24 , 24 \}_{2p} & 1+5p^{2} & \mbox{no} & 4\\ \hline
 &  &  &  &  & {\mathcal M}^- &\{ 2p , 24 \}_{24} & 2+ p^{2}( 22p- 24)/2p & & \\ \hline
5 & 4 & 2 & 24 & B & {\mathcal M}^+ & \{ 24 , 24 \}_{2p} & 1+5p^{4} & \mbox{no} & 2\\ \hline
 &  &  &  &  & {\mathcal M}^- &\{ 2p , 24 \}_{24} & 2+ p^{4}( 22p- 24)/2p & & \\ \hline
7 & 4 & 2 & 24 & B & {\mathcal M}^+ & \{ 24 , 24 \}_{2p} & 1+5p^{4} & \mbox{no} & 2\\ \hline
 &  &  &  &  & {\mathcal M}^- &\{ 2p , 24 \}_{24} & 2+ p^{4}( 22p- 24)/2p & & \\ \hline
11 & 4 & 2 & 24 & B & {\mathcal M}^+ & \{ 24 , 24 \}_{2p} & 1+5p^{4} & \mbox{yes} & 2\\ \hline
 &  &  &  &  & {\mathcal M}^- &\{ 2p , 24 \}_{24} & 2+ p^{4}( 22p- 24)/2p & & \\ \hline
13 & 4 & 2 & 24 & B & {\mathcal M}^+ & \{ 24 , 24 \}_{2p} & 1+5p^{4} & \mbox{yes} & 2\\ \hline
 &  &  &  &  & {\mathcal M}^- &\{ 2p , 24 \}_{24} & 2+ p^{4}( 22p- 24)/2p & & \\ \hline
17 & 4 & 2 & 24 & B & {\mathcal M}^+ & \{ 24 , 24 \}_{2p} & 1+5p^{4} & \mbox{no} & 2\\ \hline
 &  &  &  &  & {\mathcal M}^- &\{ 2p , 24 \}_{24} & 2+ p^{4}( 22p- 24)/2p & & \\ \hline
19 & 4 & 2 & 24 & B & {\mathcal M}^+ & \{ 24 , 24 \}_{2p} & 1+5p^{4} & \mbox{no} & 2\\ \hline
 &  &  &  &  & {\mathcal M}^- &\{ 2p , 24 \}_{24} & 2+ p^{4}( 22p- 24)/2p & & \\ \hline
23 & 2 & 1 & 24 & A & {\mathcal M}^+ & \{ 24 , 24 \}_{2p} & 1+5p^{2} & \mbox{no} & 4\\ \hline
 &  &  &  &  & {\mathcal M}^- &\{ 2p , 24 \}_{24} & 2+ p^{2}( 22p- 24)/2p & & \\ \hline
\end{array}$$

$$\begin{array}{|r|r|r|r|r|r|r|r|r|r|}\hline
p\; \mbox{mod}\; n & d & e & n & \mbox{case} & \mbox{map} & \mbox{type} & \mbox{genus}  & \mbox{self-dual} & \phi(n)/d \\ \hline
1 & 2 & 1 & 26 & B & {\mathcal M}^+ & \{ 13 , 26 \}_{2p} & 1+5p^{2} & \mbox{no} & 6\\ \hline
 &  &  &  &  & {\mathcal M}^- &\{ 2p , 26 \}_{13} & 2+ p^{2}( 24p- 26)/2p & & \\ \hline
3 & 6 & 3 & 26 & B & {\mathcal M}^+ & \{ 13 , 26 \}_{2p} & 1+5p^{6} & \mbox{no} & 2\\ \hline
 &  &  &  &  & {\mathcal M}^- &\{ 2p , 26 \}_{13} & 2+ p^{6}( 24p- 26)/2p & & \\ \hline
5 & 4 & 2 & 26 & A & {\mathcal M}^+ & \{ 13 , 26 \}_{2p} & 1+5p^{4} & \mbox{no} & 3\\ \hline
 &  &  &  &  & {\mathcal M}^- &\{ 2p , 26 \}_{13} & 2+ p^{4}( 24p- 26)/2p & & \\ \hline
7 & 12 & 6 & 26 & A & {\mathcal M}^+ & \{ 13 , 26 \}_{2p} & 1+5p^{12} & \mbox{no} & 1\\ \hline
 &  &  &  &  & {\mathcal M}^- &\{ 2p , 26 \}_{13} & 2+ p^{12}( 24p- 26)/2p & & \\ \hline
9 & 6 & 3 & 26 & B & {\mathcal M}^+ & \{ 13 , 26 \}_{2p} & 1+5p^{6} & \mbox{no} & 2\\ \hline
 &  &  &  &  & {\mathcal M}^- &\{ 2p , 26 \}_{13} & 2+ p^{6}( 24p- 26)/2p & & \\ \hline
11 & 12 & 6 & 26 & A & {\mathcal M}^+ & \{ 13 , 26 \}_{2p} & 1+5p^{12} & \mbox{no} & 1\\ \hline
 &  &  &  &  & {\mathcal M}^- &\{ 2p , 26 \}_{13} & 2+ p^{12}( 24p- 26)/2p & & \\ \hline
15 & 12 & 6 & 26 & A & {\mathcal M}^+ & \{ 13 , 26 \}_{2p} & 1+5p^{12} & \mbox{no} & 1\\ \hline
 &  &  &  &  & {\mathcal M}^- &\{ 2p , 26 \}_{13} & 2+ p^{12}( 24p- 26)/2p & & \\ \hline
17 & 6 & 3 & 26 & A & {\mathcal M}^+ & \{ 13 , 26 \}_{2p} & 1+5p^{6} & \mbox{no} & 2\\ \hline
 &  &  &  &  & {\mathcal M}^- &\{ 2p , 26 \}_{13} & 2+ p^{6}( 24p- 26)/2p & & \\ \hline
19 & 12 & 6 & 26 & A & {\mathcal M}^+ & \{ 13 , 26 \}_{2p} & 1+5p^{12} & \mbox{no} & 1\\ \hline
 &  &  &  &  & {\mathcal M}^- &\{ 2p , 26 \}_{13} & 2+ p^{12}( 24p- 26)/2p & & \\ \hline
21 & 4 & 2 & 26 & A & {\mathcal M}^+ & \{ 13 , 26 \}_{2p} & 1+5p^{4} & \mbox{no} & 3\\ \hline
 &  &  &  &  & {\mathcal M}^- &\{ 2p , 26 \}_{13} & 2+ p^{4}( 24p- 26)/2p & & \\ \hline
23 & 6 & 3 & 26 & A & {\mathcal M}^+ & \{ 13 , 26 \}_{2p} & 1+5p^{6} & \mbox{no} & 2\\ \hline
 &  &  &  &  & {\mathcal M}^- &\{ 2p , 26 \}_{13} & 2+ p^{6}( 24p- 26)/2p & & \\ \hline
25 & 2 & 1 & 26 & A & {\mathcal M}^+ & \{ 13 , 26 \}_{2p} & 1+5p^{2} & \mbox{no} & 6\\ \hline
 &  &  &  &  & {\mathcal M}^- &\{ 2p , 26 \}_{13} & 2+ p^{2}( 24p- 26)/2p & & \\ \hline
\end{array}$$

$$\begin{array}{|r|r|r|r|r|r|r|r|r|r|}\hline
p\; \mbox{mod}\; n & d & e & n & \mbox{case} & \mbox{map} & \mbox{type} & \mbox{genus}  & \mbox{self-dual} & \phi(n)/d \\ \hline
1 & 2 & 1 & 28 & B & {\mathcal M}^+ & \{ 28 , 28 \}_{2p} & 1+6p^{2} & \mbox{no} & 6\\ \hline
 &  &  &  &  & {\mathcal M}^- &\{ 2p , 28 \}_{28} & 2+ p^{2}( 26p- 28)/2p & & \\ \hline
3 & 6 & 3 & 28 & A & {\mathcal M}^+ & \{ 28 , 28 \}_{2p} & 1+6p^{6} & \mbox{no} & 2\\ \hline
 &  &  &  &  & {\mathcal M}^- &\{ 2p , 28 \}_{28} & 2+ p^{6}( 26p- 28)/2p & & \\ \hline
5 & 12 & 6 & 28 & B & {\mathcal M}^+ & \{ 28 , 28 \}_{2p} & 1+6p^{12} & \mbox{yes} & 1\\ \hline
 &  &  &  &  & {\mathcal M}^- &\{ 2p , 28 \}_{28} & 2+ p^{12}( 26p- 28)/2p & & \\ \hline
9 & 6 & 3 & 28 & B & {\mathcal M}^+ & \{ 28 , 28 \}_{2p} & 1+6p^{6} & \mbox{no} & 2\\ \hline
 &  &  &  &  & {\mathcal M}^- &\{ 2p , 28 \}_{28} & 2+ p^{6}( 26p- 28)/2p & & \\ \hline
11 & 12 & 6 & 28 & B & {\mathcal M}^+ & \{ 28 , 28 \}_{2p} & 1+6p^{12} & \mbox{yes} & 1\\ \hline
 &  &  &  &  & {\mathcal M}^- &\{ 2p , 28 \}_{28} & 2+ p^{12}( 26p- 28)/2p & & \\ \hline
13 & 4 & 2 & 28 & B & {\mathcal M}^+ & \{ 28 , 28 \}_{2p} & 1+6p^{4} & \mbox{yes} & 3\\ \hline
 &  &  &  &  & {\mathcal M}^- &\{ 2p , 28 \}_{28} & 2+ p^{4}( 26p- 28)/2p & & \\ \hline
15 & 4 & 2 & 28 & B & {\mathcal M}^+ & \{ 28 , 28 \}_{2p} & 1+6p^{4} & \mbox{yes} & 3\\ \hline
 &  &  &  &  & {\mathcal M}^- &\{ 2p , 28 \}_{28} & 2+ p^{4}( 26p- 28)/2p & & \\ \hline
17 & 12 & 6 & 28 & B & {\mathcal M}^+ & \{ 28 , 28 \}_{2p} & 1+6p^{12} & \mbox{yes} & 1\\ \hline
 &  &  &  &  & {\mathcal M}^- &\{ 2p , 28 \}_{28} & 2+ p^{12}( 26p- 28)/2p & & \\ \hline
19 & 6 & 3 & 28 & A & {\mathcal M}^+ & \{ 28 , 28 \}_{2p} & 1+6p^{6} & \mbox{no} & 2\\ \hline
 &  &  &  &  & {\mathcal M}^- &\{ 2p , 28 \}_{28} & 2+ p^{6}( 26p- 28)/2p & & \\ \hline
23 & 12 & 6 & 28 & B & {\mathcal M}^+ & \{ 28 , 28 \}_{2p} & 1+6p^{12} & \mbox{yes} & 1\\ \hline
 &  &  &  &  & {\mathcal M}^- &\{ 2p , 28 \}_{28} & 2+ p^{12}( 26p- 28)/2p & & \\ \hline
25 & 6 & 3 & 28 & B & {\mathcal M}^+ & \{ 28 , 28 \}_{2p} & 1+6p^{6} & \mbox{no} & 2\\ \hline
 &  &  &  &  & {\mathcal M}^- &\{ 2p , 28 \}_{28} & 2+ p^{6}( 26p- 28)/2p & & \\ \hline
27 & 2 & 1 & 28 & A & {\mathcal M}^+ & \{ 28 , 28 \}_{2p} & 1+6p^{2} & \mbox{no} & 6\\ \hline
 &  &  &  &  & {\mathcal M}^- &\{ 2p , 28 \}_{28} & 2+ p^{2}( 26p- 28)/2p & & \\ \hline
\end{array}$$

$$\begin{array}{|r|r|r|r|r|r|r|r|r|r|}\hline
p\; \mbox{mod}\; n & d & e & n & \mbox{case} & \mbox{map} & \mbox{type} & \mbox{genus}  & \mbox{self-dual} & \phi(n)/d \\ \hline
1 & 2 & 1 & 30 & B & {\mathcal M}^+ & \{ 15 , 30 \}_{2p} & 1+6p^2 & \mbox{no} & 4\\ \hline
 &  &  &  &  & {\mathcal M}^- &\{ 2p , 30 \}_{15} & 2+ p^2( 28p- 30)/2p & & \\ \hline
7 & 8 & 4 & 30 & B & {\mathcal M}^+ & \{ 15 , 30 \}_{2p} & 1+6p^8 & \mbox{no} & 1\\ \hline
 &  &  &  &  & {\mathcal M}^- &\{ 2p , 30 \}_{15} & 2+ p^8( 28p- 30)/2p & & \\ \hline
11 & 4 & 2 & 30 & B & {\mathcal M}^+ & \{ 15 , 30 \}_{2p} & 1+6p^4 & \mbox{no} & 2\\ \hline
 &  &  &  &  & {\mathcal M}^- &\{ 2p , 30 \}_{15} & 2+ p^4( 28p- 30)/2p & & \\ \hline
13 & 8 & 4 & 30 & B & {\mathcal M}^+ & \{ 15 , 30 \}_{2p} & 1+6p^8 & \mbox{no} & 1\\ \hline
 &  &  &  &  & {\mathcal M}^- &\{ 2p , 30 \}_{15} & 2+ p^8( 28p- 30)/2p & & \\ \hline
17 & 8 & 4 & 30 & B & {\mathcal M}^+ & \{ 15 , 30 \}_{2p} & 1+6p^8 & \mbox{no} & 1\\ \hline
 &  &  &  &  & {\mathcal M}^- &\{ 2p , 30 \}_{15} & 2+ p^8( 28p- 30)/2p & & \\ \hline
19 & 4 & 2 & 30 & B & {\mathcal M}^+ & \{ 15 , 30 \}_{2p} & 1+6p^4 & \mbox{no} & 2\\ \hline
 &  &  &  &  & {\mathcal M}^- &\{ 2p , 30 \}_{15} & 2+ p^4( 28p- 30)/2p & & \\ \hline
23 & 8 & 4 & 30 & B & {\mathcal M}^+ & \{ 15 , 30 \}_{2p} & 1+6p^8 & \mbox{no} & 1\\ \hline
 &  &  &  &  & {\mathcal M}^- &\{ 2p , 30 \}_{15} & 2+ p^8( 28p- 30)/2p & & \\ \hline
29 & 2 & 1 & 30 & A & {\mathcal M}^+ & \{ 15 , 30 \}_{2p} & 1+6p^2 & \mbox{no} & 4\\ \hline
 &  &  &  &  & {\mathcal M}^- &\{ 2p , 30 \}_{15} & 2+ p^2( 28p- 30)/2p & & \\ \hline
\end{array}$$

\newpage

The following tables give all information for each pair $p$ and $n$ in one row, with the type and genus of ${\mathcal M}^+$ followed by those for its Petrie dual ${\mathcal M}^-$.

$$\begin{array}{|r|r|r|r|r|r|r|r|r|r|r|}\hline
p& d & e & n & \mbox{case} & \mbox{type} & \mbox{genus}  & \mbox{self dual} & \mbox{type} & \mbox{genus}  & \phi(n)/d \\ \hline
2 & 2 & 1 & 3 & A & \{ 3 , 3 \}_{4} & 0 & \mbox{yes} & \{ 4 , 3 \}_{3} & 1 & 1\\ \hline
2 & 4 & 2 & 5 & A & \{ 5 , 5 \}_{4} & 5 & \mbox{yes} & \{ 4 , 5 \}_{5} & 6 & 1\\ \hline
2 & 6 & 3 & 7 & B & \{ 7 , 7 \}_{4} & 49 & \mbox{yes} & \{ 4 , 7 \}_{7} & 50 & 1\\ \hline
2 & 6 & 3 & 9 & A & \{ 9 , 9 \}_{4} & 81 & \mbox{yes} & \{ 4 , 9 \}_{9} & 82 & 1\\ \hline
2 & 10 & 5 & 11 & A & \{ 11 , 11 \}_{4} & 1793 & \mbox{yes} & \{ 4 , 11 \}_{11} & 1794 & 1\\ \hline
2 & 12 & 6 & 13 & A & \{ 13 , 13 \}_{4} & 9217 & \mbox{yes} & \{ 4 , 13 \}_{13} & 9218 & 1\\ \hline
2 & 8 & 4 & 15 & B & \{ 15 , 15 \}_{4} & 705 & \mbox{yes} & \{ 4 , 15 \}_{15} & 706 & 1\\ \hline
2 & 8 & 4 & 17 & A & \{ 17 , 17 \}_{4} & 833 & \mbox{yes} & \{ 4 , 17 \}_{17} & 834 & 2\\ \hline
2 & 18 & 9 & 19 & A & \{ 19 , 19 \}_{4} & 983041 & \mbox{yes} & \{ 4 , 19 \}_{19} & 983042 & 1\\ \hline
2 & 12 & 6 & 21 & B & \{ 21 , 21 \}_{4} & 17409 & \mbox{yes} & \{ 4 , 21 \}_{21} & 17410 & 1\\ \hline
2 & 22 & 11 & 23 & B & \{ 23 , 23 \}_{4} & 19922945 & \mbox{yes} & \{ 4 , 23 \}_{23} & 19922946 & 1\\ \hline
2 & 20 & 10 & 25 & A & \{ 25 , 25 \}_{4} & 5505025 & \mbox{yes} & \{ 4 , 25 \}_{25} & 5505026 & 1\\ \hline
2 & 18 & 9 & 27 & A & \{ 27 , 27 \}_{4} & 1507329 & \mbox{yes} & \{ 4 , 27 \}_{27} & 1507330 & 1\\ \hline
2 & 28 & 14 & 29 & A & \{ 29 , 29 \}_{4} & 1677721601 & \mbox{yes} & \{ 4 , 29 \}_{29} & 1677721602 & 1\\ \hline
2 & 10 & 5 & 31 & B & \{ 31 , 31 \}_{4} & 6913 & \mbox{yes} & \{ 4 , 31 \}_{31} & 6914 & 3\\ \hline
2 & 10 & 5 & 33 & A & \{ 33 , 33 \}_{4} & 7425 & \mbox{yes} & \{ 4 , 33 \}_{33} & 7426 & 2\\ \hline
2 & 24 & 12 & 35 & B & \{ 35 , 35 \}_{4} & 130023425 & \mbox{yes} & \{ 4 , 35 \}_{35} & 130023426 & 1\\ \hline
2 & 36 & 18 & 37 & A & \{ 37 , 37 \}_{4} & 566935683073 & \mbox{yes} & \{ 4 , 37 \}_{37} & 566935683074 & 1\\ \hline
2 & 24 & 12 & 39 & B & \{ 39 , 39 \}_{4} & 146800641 & \mbox{yes} & \{ 4 , 39 \}_{39} & 146800642 & 1\\ \hline
2 & 20 & 10 & 41 & A & \{ 41 , 41 \}_{4} & 9699329 & \mbox{yes} & \{ 4 , 41 \}_{41} & 9699330 & 2\\ \hline
2 & 14 & 7 & 43 & A & \{ 43 , 43 \}_{4} & 159745 & \mbox{yes} & \{ 4 , 43 \}_{43} & 159746 & 3\\ \hline
2 & 24 & 12 & 45 & B & \{ 45 , 45 \}_{4} & 171966465 & \mbox{yes} & \{ 4 , 45 \}_{45} & 171966466 & 1\\ \hline
2 & 46 & 23 & 47 & B & \{ 47 , 47 \}_{4} & 756463999909889 & \mbox{yes} & \{ 4 , 47 \}_{47} & 756463999909890 & 1\\ \hline
2 & 42 & 21 & 49 & B & \{ 49 , 49 \}_{4} & 49478023249921 &\mbox{yes} & \{ 4 , 49 \}_{49} & 49478023249922 & 1\\ \hline
\end{array}$$

$$\begin{array}{|r|r|r|r|r|r|r|r|r|r|r|}\hline
p\; \mbox{mod}\; n & d & e & n & \mbox{case} & \mbox{type} & \mbox{genus}  & \mbox{self dual} & \mbox{type} & \mbox{genus}  & \phi(n)/d \\ \hline
1 & 2 & 1 & 4 & B & \{ 4 , 4 \}_{2p} & 1 & \mbox{yes} & \{ 2p , 4 \}_{4} & 2+ p^2( 2p- 4)/2p & 1\\ \hline
3 & 2 & 1 & 4 & A & \{ 4 , 4 \}_{2p} & 1 & \mbox{yes} & \{ 2p , 4 \}_{4} & 2+ p^2( 2p- 4)/2p & 1\\ \hline
\end{array}$$

$$\begin{array}{|r|r|r|r|r|r|r|r|r|r|r|}\hline
p\; \mbox{mod}\; n & d & e & n & \mbox{case} & \mbox{type} & \mbox{genus}  & \mbox{self dual} & \mbox{type} & \mbox{genus}  & \phi(n)/d \\ \hline
1 & 2 & 1 & 28 & B & \{ 28 , 28 \}_{2p} & 1+6p^2 & \mbox{no} & \{ 2p , 28 \}_{28} & 2+ p^2( 26p- 28)/2p & 6\\ \hline
3 & 6 & 3 & 28 & A & \{ 28 , 28 \}_{2p} & 1+6p^6 & \mbox{no} & \{ 2p , 28 \}_{28} & 2+ p^6( 26p- 28)/2p & 2\\ \hline
5 & 12 & 6 & 28 & B & \{ 28 , 28 \}_{2p} & 1+6p^{12} & \mbox{yes} & \{ 2p , 28 \}_{28} & 2+ p^12( 26p- 28)/2p & 1\\ \hline
9 & 6 & 3 & 28 & B & \{ 28 , 28 \}_{2p} & 1+6p^6 & \mbox{no} & \{ 2p , 28 \}_{28} & 2+ p^6( 26p- 28)/2p & 2\\ \hline
11 & 12 & 6 & 28 & B & \{ 28 , 28 \}_{2p} & 1+6p^{12} & \mbox{yes} & \{ 2p , 28 \}_{28} & 2+ p^12( 26p- 28)/2p & 1\\ \hline
13 & 4 & 2 & 28 & B & \{ 28 , 28 \}_{2p} & 1+6p^4 & \mbox{yes} & \{ 2p , 28 \}_{28} & 2+ p^4( 26p- 28)/2p & 3\\ \hline
15 & 4 & 2 & 28 & B & \{ 28 , 28 \}_{2p} & 1+6p^4 & \mbox{yes} & \{ 2p , 28 \}_{28} & 2+ p^4( 26p- 28)/2p & 3\\ \hline
17 & 12 & 6 & 28 & B & \{ 28 , 28 \}_{2p} & 1+6p^{12} & \mbox{yes} & \{ 2p , 28 \}_{28} & 2+ p^12( 26p- 28)/2p & 1\\ \hline
19 & 6 & 3 & 28 & A & \{ 28 , 28 \}_{2p} & 1+6p^6 & \mbox{no} & \{ 2p , 28 \}_{28} & 2+ p^6( 26p- 28)/2p & 2\\ \hline
23 & 12 & 6 & 28 & B & \{ 28 , 28 \}_{2p} & 1+6p^{12} & \mbox{yes} & \{ 2p , 28 \}_{28} & 2+ p^12( 26p- 28)/2p & 1\\ \hline
25 & 6 & 3 & 28 & B & \{ 28 , 28 \}_{2p} & 1+6p^6 & \mbox{no} & \{ 2p , 28 \}_{28} & 2+ p^6( 26p- 28)/2p & 2\\ \hline
27 & 2 & 1 & 28 & A & \{ 28 , 28 \}_{2p} & 1+6p^2 & \mbox{no} & \{ 2p , 28 \}_{28} & 2+ p^2( 26p- 28)/2p & 6\\ \hline
\end{array}$$

$$\begin{array}{|r|r|r|r|r|r|r|r|r|r|r|}\hline
p\; \mbox{mod}\; n & d & e & n & \mbox{case} & \mbox{type} & \mbox{genus}  & \mbox{self dual} & \mbox{type} & \mbox{genus}  & \phi(n)/d \\ \hline
1 & 2 & 1 & 30 & B & \{ 15 , 30 \}_{2p} & 1+6p^2 & \mbox{no} & \{ 2p , 30 \}_{15} & 2+ p^2( 28p- 30)/2p & 4\\ \hline
7 & 8 & 4 & 30 & B & \{ 15 , 30 \}_{2p} & 1+6p^8 & \mbox{no} & \{ 2p , 30 \}_{15} & 2+ p^8( 28p- 30)/2p & 1\\ \hline
11 & 4 & 2 & 30 & B & \{ 15 , 30 \}_{2p} & 1+6p^4 & \mbox{no} & \{ 2p , 30 \}_{15} & 2+ p^4( 28p- 30)/2p & 2\\ \hline
13 & 8 & 4 & 30 & B & \{ 15 , 30 \}_{2p} & 1+6p^8 & \mbox{no} & \{ 2p , 30 \}_{15} & 2+ p^8( 28p- 30)/2p & 1\\ \hline
17 & 8 & 4 & 30 & B & \{ 15 , 30 \}_{2p} & 1+6p^8 & \mbox{no} & \{ 2p , 30 \}_{15} & 2+ p^8( 28p- 30)/2p & 1\\ \hline
19 & 4 & 2 & 30 & B & \{ 15 , 30 \}_{2p} & 1+6p^4 & \mbox{no} & \{ 2p , 30 \}_{15} & 2+ p^4( 28p- 30)/2p & 2\\ \hline
23 & 8 & 4 & 30 & B & \{ 15 , 30 \}_{2p} & 1+6p^8 & \mbox{no} & \{ 2p , 30 \}_{15} & 2+ p^8( 28p- 30)/2p & 1\\ \hline
29 & 2 & 1 & 30 & A & \{ 15 , 30 \}_{2p} & 1+6p^2 & \mbox{no} & \{ 2p , 30 \}_{15} & 2+ p^2( 28p- 30)/2p & 4\\ \hline
\end{array}$$

%\fi

\end{document}